\documentclass[10pt,reqno]{amsart}

\usepackage{fixltx2e}                       
\usepackage[usenames,dvipsnames]{xcolor}
\usepackage{fancyhdr}
\usepackage{amsmath,amsfonts,amsbsy,amsgen,amscd,amssymb,amsthm}
\usepackage{mathrsfs,eucal}
\usepackage{mathtools}
\usepackage[font=small,margin=0.25in,labelfont={sc},labelsep={colon}]{caption}
\usepackage{tikz}
\usepackage{enumitem}
\usepackage[colorlinks=true,linkcolor=blue,citecolor=green]{hyperref}
\usepackage{graphicx}

\usepackage[T1]{fontenc}
\usepackage{bm}

\newtheorem{theorem}{Theorem}[section]
\newtheorem{lemma}[theorem]{Lemma}
\newtheorem{corollary}[theorem]{Corollary}
\newtheorem{definition}[theorem]{Definition}
\newtheorem{proposition}[theorem]{Proposition}
\newtheorem{remark}[theorem]{Remark}
\newtheorem*{theorem*}{Informal Statement of Theorem}

\DeclareMathOperator{\eps}{\varepsilon}
\DeclareMathOperator{\EE}{\mathbb{E}}
\DeclareMathOperator{\PP}{\mathbb{P}}
\DeclareMathOperator{\RR}{\mathbb{R}}

\DeclareMathOperator{\III}{\mathcal{I}}
\DeclareMathOperator{\JJJ}{\mathcal{J}}
\DeclareMathOperator{\GGG}{\mathcal{G}}
\DeclareMathOperator{\DDD}{\mathcal{D}}
\DeclareMathOperator{\EEE}{\mathcal{E}}
\DeclareMathOperator{\ZZ}{\mathbb{Z}}
\DeclareMathOperator{\Tr}{\mathrm{Tr}}
\DeclareMathOperator{\1}{\mathbf{1}}

\newcommand{\diag}{\mathrm{diag}}
\newcommand{\rank}{\mathrm{rank}}

\begin{document}

\title[Random Laplacian matrices and convex relaxations]{Random Laplacian matrices and convex relaxations}
\author{Afonso S. Bandeira}
\address{Department of Mathematics, Massachusetts Institute of Technology, Cambridge, MA 02139, USA}
\thanks{This work was carried out while the author was with the Program in Applied and Computational Mathematics in Princeton 
University and supported by AFOSR Grant No.\ FA9550-12-1-0317.}
\email{bandeira@mit.edu}
\date{\today}

\begin{abstract}
The largest eigenvalue of a matrix is always larger or equal than its largest diagonal entry. We show that for a large class of random Laplacian matrices, this bound is essentially tight: the largest eigenvalue is, up to lower order terms, often the size of the largest diagonal entry.

Besides being a simple tool to obtain precise estimates on the largest eigenvalue of a large class of random Laplacian matrices, our main result settles a number of open problems related to the tightness of certain convex relaxation--based algorithms. It easily implies the optimality of the semidefinite relaxation approaches to problems such as $\ZZ_2$ Synchronization and Stochastic Block Model recovery. Interestingly, this result readily implies the connectivity threshold for Erd\H{o}s--R\'{e}nyi graphs and suggests that these three phenomena are manifestations of the same underlying principle.
The main tool is a recent estimate on the spectral norm of matrices with independent entries by van~Handel and the author.
\end{abstract}

%

\maketitle

\thispagestyle{empty}


\section{Introduction}

Towards the end of the 1950s, Eugene Wigner~\cite{Wigner_1958} made the remarkable finding that the spectrum of a large class of random matrices is, in high dimension, distributed essentially the same way: under mild assumptions, the distribution of the spectrum converges to the so-called \emph{Wigner semicircle law}. The study of spectral properties of random matrices has since spawned a panoply of fascinating research with important implications in many areas. We refer the reader to the books~\cite{Tao_topicsRMT,Anderson_Guionnet_Zeitouni_IntroRandomMatrices} for more on this subject.

The present paper addresses the problem of estimating the largest eigenvalue of a large class of Laplacian matrices. The investigation of such problems has strong motivations from algorithmic analysis. Indeed, the performance of many popular algorithms is tightly connected with the largest eigenvalue of some matrix that depends on its input, and so studying the performance of such algorithms over random inputs involves understanding the behavior of the largest eigenvalue of a random matrix. In fact, as we will see, the estimates derived here play a crucial role in understanding the typical performance of a natural semidefinite programming--based approach for solving certain computationally hard problems on graphs, such as community detection.

We use the term \emph{Laplacian matrix} to refer to symmetric matrices whose rows and columns sum to zero. While oftentimes Laplacians are also thought of as being positive semidefinite, the matrices we will treat will not necessarily satisfy this property. Spectral graph theory inspires a useful way of thinking about these matrices~\cite{FanChung_SpectralGraphTheory}. Given a graph on $n$ nodes with edge set $E$, its adjacency matrix $A\in\RR^{n\times n}$ is defined by $A_{ij}=1$ if $(i,j)\in E$ and $A_{ij}=0$ otherwise, and its degree matrix $D_A$ is a diagonal matrix whose $i$-th diagonal entry is equal to the degree of node $i$. The Laplacian of the graph is defined to be $L_A = D_A - A$. The spectrum of the graph Laplacian matrix is known to contain important information about the graph~\cite{FanChung_SpectralGraphTheory}, and has been studied for random graphs~\cite{Ding_RandomLaplacians,Chung_ComplexGraphsNetworks,Bryc_06_MarkovMatrices}. Analogously, we make the following definition.

\begin{definition}\label{def:Laplacian_and_Degree}
Given a symmetric matrix $X\in \RR^{n\times n}$, we define the Laplacian $L_X$ of $X$ as $$L_X = D_X - X,$$ where $D_X$ is the diagonal matrix whose diagonal entries are given by
\[
\left(D_X\right)_{ii} = \sum_{j=1}^n X_{ij}.
\]
We will refer to any such matrix $L_X$ as a Laplacian matrix. Note that these are precisely the symmetric matrices $L$ for which $L\1$ = 0, where $\1\in\RR^n$ denotes the all-ones vector.
\end{definition}

This paper is concerned with a class of random Laplacian matrices $L_X$ where the entries of the matrix $X$ are independent centered (but not necessarily identically distributed) random variables. Our main result is that, under mild and easily verifiable conditions,  the largest eigenvalue of $L_X$ is, up to lower order terms, given by its largest diagonal entry. While we defer the formal statement of our main results~\footnote{Our results will be of nonasymptotic nature (we refer the interested reader to~\cite{VershyninNARandomMatrices} for a tutorial on nonasymptotic estimates in random matrix theory).} to Section~\ref{section:mainresults}, we informaly state them here.

\begin{theorem*}[\ref{theorem:main}]
 
Let $L$ be an $n\times n$ symmetric random Laplacian matrix (i.e. satisfying $L \1=0$) with centered independent off-diagonal entries such that
\(\sum_{j\in [n]\setminus i}\EE L_{ij}^2 \) is equal for every $i$, and 
\[
 \sum_{j\in [n]\setminus i}\EE L_{ij}^2 \gtrsim  \max_{i\neq j}\left\| L_{ij} \right\|_{\infty}^2 \log n.
\]
Then, with high probability,
\[
 \lambda_{\max}(L) - \max_{i}L_{ii} \lesssim (\log n)^{-\frac12} \max_{i}L_{ii}.
\]

\end{theorem*}

Not only does our main result provide an extremely simple tool to precisely estimate the largest eigenvalue of Laplacian matrices, but in the applications studied below, the largest diagonal value also enjoys an interpretation that is intimately tied to the underlying problem.

To illustrate the latter point, we turn back to graph theory. It is well known that the spectrum of the Laplacian of a graph dictates whether or not the graph is connected. On the other hand, its diagonal is simply given by the degrees of the nodes of the graph. A relation between the spectrum of the Laplacian and its diagonal could then translate into a relation between degrees of nodes of a graph and its connectivity. In fact, such a relation is already known to exist: \textit{The phase transition for connectivity of Erd\H{o}s--R\'{e}nyi graphs\footnote{The Erd\H{o}s--R\'{e}nyi model for random graphs will be discussed in more detail in Section~\ref{subsection:ErdosRenyi}.} coincides with the one for the existence of isolated nodes}. While it is true that any graph with an isolated node (a node with degree zero) cannot be connected, the converse is far from true, rendering this phenomenon particularly interesting. In Section~\ref{subsection:ErdosRenyi}, we will use our main result to provide a simple and illustrative proof of this phenomenon.

We will use our main result to give sharp guarantees for certain algorithms that solve the $\ZZ_2$ Synchronization problem and the community detection problem in the Stochastic Block Model. The $\ZZ_2$ Synchronization problem consists of recovering binary labels $x_i=\pm 1$ associated with nodes of a graph from noisy (pairwise) measurements of $x_ix_j$ whenever $(i,j)$ is an edge of the graph (see~\cite{ASinger_2011_angsync}). This problem is intimately related to correlation clustering~\cite{correlation_clustering_04}. Despite its hardness, spectral methods and semidefinite programming--based methods are known to perform well in both the worst-case~\cite{Bandeira_Singer_Spielman_OdCheeger} and average-case settings~\cite{Abbe_Z2Synch,Abbe_Z2SynchER,Cucuringu_Z2Synch}.\footnote{The information-theoretic limits of this problem have also been investigated~\cite{Abbe_Z2Synch,Abbe_Z2SynchER,Chen_Goldsmith_ISIT2014,Chen_Goldsmith_2015}.}

Community detection, or clustering, in a graph is a central problem in countless applications. Unfortunately, even the simplified version of partitioning a graph into two vertex sets, with the same size, that minimize the number of edges across the partition, referred to as \emph{minimum bisection}, is known to be NP-hard. Nevertheless, certain heuristics are known to work well for typical realizations of random graph models that exhibit community structure~\cite{McSherry_SBM1,Boppana_SBM1,Feige_Kilian_bisection_01}. In this setting, a particularly popular model is the Stochastic Block Model with two communities.

\begin{definition}\label{def:SBM}(Stochastic Block Model with two communities)
 Given $n$ even, and $0\leq p,q \leq 1$, we say that a random graph $G$ is drawn from $\GGG(n,p,q)$, the Stochastic Block Model with two communities, if $G$ has $n$ nodes, divided in two clusters of $\frac{n}2$ nodes each, and for each pair of vertices $i,j$, $(i,j)$ is an edge of $G$ with probability $p$ if $i$ and $j$ are in the same cluster and with probability $q$ otherwise, independently from any other edge. 
\end{definition}

We will focus on the setting $p>q$. The problem of recovering, from a realization $G\sim\GGG(n,p,q)$, the original partition of the underlying vertices gained popularity when Decelle et al.~\cite{Decelle_SBM} conjectured a fascinating phase transition in the constant average--degree regime. More precisely, if $p = \frac{a}n$ and $q = \frac{b}n$ with $a>b$ constants, it was conjectured that as long as 
\[
 (a-b)^2 > 2(a+b),
\]
it is possible to make an estimate of the original partition that correlates with the true partition, and that below this threshold it is impossible to do so. This conjecture was later proven in a remarkable series of works by Mossel et al.~\cite{Mossel_SBM1,Mossel_SBM2} and Massoulie~\cite{Massoulie_SBM}. Instead of settling for an estimate that correlates with the true partition, we will focus on exactly recovering the partition. A phase transition for this problem was established by Abbe et al.~\cite{SBM_paper_arxiv_v4} and independently by Mossel et al.~\cite{Mossel_SBM3_exact}. We will show that a certain semidefinite programming based--algorithm succeeds up to the information theoretical threshold, thus settling a problem posed in~\cite{SBM_paper_arxiv_v4}. We remark that, while the present paper was being written, it was brought to our attention that this problem was also solved independently by parallel research efforts of Hajek et al.~\cite{Hajek_et_al_SBM_SDP}.

The use of semidefinite relaxations in combinatorial optimization dates back to the late 1970s with the seminal work of L\'aszl\'o Lov\'asz~\cite{Lovasz_ShannonCapacity} in the so-called \emph{Lov\'asz theta function}, this approach was shortly after made algorithmic in~\cite{Grotschel_Lovasz_Schrijver_ellipsoidmethod}. In the first half of the 1990s, interior point methods were adapted to solve semidefinite programs~\cite{Alizadeh_93,Nesterov_Nemirovskii_interiorpointmethods_94}, providing reasonably efficient methods to solve these problems. In 1995, Goemans and Williamson, devised the first approximation algorithm based on semidefinite programming~\cite{MXGoemans_DPWilliamson_1995}. Their algorithm gave the best known approximation ratio to the \texttt{Max-Cut} problem. Ever since, many approximation algorithms have been designed based on semidefinite programming. In fact, the algorithm we will analyze is greatly inspired by the semidefinite relaxation in~\cite{MXGoemans_DPWilliamson_1995}. Remarkably, an important conjecture of Khot~\cite{SKhot_2002} is known to imply that for a large class of problems including \texttt{Max-Cut}, this approach produces optimal approximation ratios~\cite{Raghavendra_2008_optimalitySDP_UG}.

An approximation ratio is a guarantee that, for any possible instance of the input, the algorithm outputs a solution whose performance is at least a certain fraction (the approximation ratio) of the optimal one. The worst-case nature of this type of guarantee is often pessimistic. A popular alternative is to equip the input with a distribution (such as, for example, the Stochastic Block Model) and give guarantees for most inputs. More precisely, we will be interested in understanding when is it the case that the semidefinite relaxation approach gives exactly the correct answer (for most inputs). The tendency for a large class of semidefinite relaxations to be tight\footnote{When the optimal solution of a semidefinite relaxation is the optimal solution of the original problem we say that the relaxation is \emph{tight}.} has been observed and conjectured, for example, in~\cite{Bandeira2014open}. One of the main insights of this paper is the fact that the phenomenon described by our main result provides a unifying principle for understanding the tightness of many convex relaxations.

\subsection{Notation}

We will make use of several standard matrix and probability notations. For $M$ a matrix we will denote its $k$-th smallest eigenvalue by $\lambda_{k}(M)$, largest eigenvalue by $\lambda_{\max}(M)$, and its spectral norm by $\|M\|$. $\diag(M)$ will be used to refer to a vector with the diagonal elements of $M$ as entries. For $x\in\RR^n$ a vector, $\diag(x)$ will denote a diagonal matrix $D\in\RR^{n\times n}$ with $D_{ii}=x_i$.

$\1$ will denote the all-ones vector, whenever there is no risk of ambiguity for its dimension.

For a scalar random variable $Y$, we will write its $p$-norm as $\|Y\|_p = \left( \EE |Y|^p \right)^{1/p}$ and infinity norm as $\|Y\|_{\infty} = \inf\left\{a :\  |Y| \leq a \text{ a.\ s.} \right\}$.

Given a graph, $\deg(i)$ will be used to denote the degree of node $i$. In the case of the Stochastic Block Model, $\deg_{in}(i)$ will be used for inner-cluster degree and $\deg_{out}(i)$ for outer-cluster degree.

We will say that an event $\EEE$ happens with high probability when
\[
 \PP\left[ \EEE \right] = 1-n^{-\Omega(1)},
\]
where $n$ is an underlying parameter that is thought of going to infinity (such as the dimension of the matrices or the number of nodes in the graphs being studied).

\section{A simpler problem: \(\mathbb{Z}_2\) Synchronization with gaussian noise}\label{section:firstexample}

Before presenting our main results in Section~\ref{section:mainresults}, we will motivate them through a simplified version of the problems of $\ZZ_2$ Synchronization and recovery in the Stochastic Block Model: given a noise level $\sigma$ and a vector $z\in\{\pm1\}^n$ suppose we are given noisy measurements
\[
 Y_{ij} = z_iz_j + \sigma W_{ij},
\]
for each pair $(i,j)$, where $W_{ij}$ are i.i.d. standard gaussian random variables (with $W_{ij}=W_{ji}$). A version of this problem, over the complex numbers, is treated in~\cite{Bandeira_rankrecoveryangsynch}. Our objective is to devise an algorithm that recovers the correct $z$ with high probability. By definition, the maximum a posteriori (MAP) estimator maximizes the probability of recovering the correct variable $z$. Given that we have no a priori information on $z$ we assume a uniform prior, in that case the MAP estimator coincides with the Maximum Likelihood Estimator (MLE) for $z$. The latter is the solution of
\begin{equation}\label{littleGrothendieck:Z2:1}
\begin{array}{cl}
 \max & x^TYx \\
 \text{ s.t. } & x\in\RR^n \\
		& x_i^2 = 1,
\end{array}
\end{equation}
which is referred to as the little Grothendieck problem over $\RR$ and known to be NP-hard in general. In fact, \eqref{littleGrothendieck:Z2:1} includes the \texttt{Max-Cut} problem by taking $Y$ to be the Laplacian of a graph. In the spirit of the relaxation proposed in~\cite{MXGoemans_DPWilliamson_1995} for the \texttt{Max-Cut} problem, we take $X = xx^T$ and rewrite \eqref{littleGrothendieck:Z2:1} as
\begin{equation}\label{littleGrothendieck:Z2:1lifted}
\begin{array}{cl}
 \max & \Tr(YX) \\
 \text{ s.t. }  & X_{ii} = 1\\
		& X\succeq 0 \\
		& \rank(X) = 1.
\end{array}
\end{equation}

We now relax the nonconvex rank constraint and arrive at the following semidefinite program, which can be solved in polynomial time up to arbitrary precision~\cite{LVanderberghe_SBoyd_1996}.
\begin{equation}\label{eq:thegeneralSDPformulation}
\begin{array}{cl}
 \max & \Tr(YX) \\
 \text{ s.t. }  & X_{ii} = 1\\
		& X\succeq 0.
\end{array}
\end{equation}

As it will be clear in the proceeding sections, this relaxation is also used to solve $\ZZ_2$ Synchronization and recovery in the Stochastic Block Model, albeit for a different coefficient matrix $Y$.

In what follows we will derive conditions for when a certain rank $1$ matrix is the unique optimal solution of \eqref{eq:thegeneralSDPformulation}. Note that if $X=xx^T$ is the unique solution to \eqref{eq:thegeneralSDPformulation}, then $x$ must be the solution to \eqref{littleGrothendieck:Z2:1}, meaning that we are able to compute the MLE efficiently by solving \eqref{eq:thegeneralSDPformulation}. This motivates us to understand when is it the case that $X = xx^T$ is the unique optimal solution of \eqref{eq:thegeneralSDPformulation}. A fruitful way of approaching this relies on duality. The dual of \eqref{eq:thegeneralSDPformulation} is given by:
\begin{equation}\label{eq:thegeneralSDPformulation_DUAL}
\begin{array}{cl}
 \min & \Tr(D) \\
 \text{ s.t. }  & D \text{ is diagonal}\\
		& D - Y\succeq 0.
\end{array}
\end{equation}

Weak duality guarantees that if $X$ and $D$ are feasible solutions of respectively \eqref{eq:thegeneralSDPformulation} and \eqref{eq:thegeneralSDPformulation_DUAL} then $\Tr(YX)\leq \Tr(D)$. Indeed, since $X$ and $D-Y$ are both positive semidefinite, we must have
\begin{equation}\label{eq:weakduality:CompSlack}
 0 \leq \Tr\left[ (D-Y)X \right] = \Tr(D) - \Tr(YX).
\end{equation}
This means that if we are able to find a so-called dual certificate, a matrix $D$ feasible for \eqref{eq:thegeneralSDPformulation_DUAL} for which $\Tr(D) = \Tr(Yxx^T)$, then it guarantees that $X=xx^T$ is an optimal solution of \eqref{eq:thegeneralSDPformulation}. To guarantee uniqueness it suffices to further ensure that $\lambda_2(D-Y)>0$. In fact, if there existed another optimal solution $X$, by \eqref{eq:weakduality:CompSlack}, one would have $\Tr\left[ (D-Y)X \right]=0$ which can be shown to imply (see, for example,~\cite{Abbe_Z2Synch}), together with the feasibility of $X$, that $X=xx^T$. This establishes the following Lemma.

\begin{lemma}\label{lemma:dualcertificate}[Dual Certificate]
 Let $Y$ be a symmetric $n\times n$ matrix and $x\in\{\pm1\}^n$. If there exists a diagonal matrix $D$, such that $\Tr(D) = x^TYx$, $D-Y\succeq 0$, and $\lambda_2(D-Y)>0$ then $X=xx^T$ is the unique optimal solution of \eqref{eq:thegeneralSDPformulation}.
\end{lemma}

We take a candidate dual certificate $D$ whose diagonal elements are given by
\[
 D_{ii} = \sum_{j=1}^n Y_{ij}x_ix_j.
\]
Note that $D = D_{[\diag(x)Y\diag(x)]}$ as per Definition~\ref{def:Laplacian_and_Degree}. It is easy to see that $\Tr(D) = x^TYx$ and $(D-Y)x = 0$ which gives the following Lemma.

\begin{lemma}\label{lemma:dualcertificate_constructed}
 Let $Y$ be a symmetric $n\times n$ matrix and $x\in\{\pm1\}^n$. Let $D$ be the diagonal matrix defined as $D =  D_{[\diag(x)Y\diag(x)]}$. As long as 
\[
 \lambda_2(D-Y)>0,
\]
$X=xx^T$ is the unique optimal solution of \eqref{eq:thegeneralSDPformulation}.
\end{lemma}

Note that these guarantees, (Lemmas~\ref{lemma:dualcertificate} and~\ref{lemma:dualcertificate_constructed}) do not depend on the matrix $Y$ or the distribution from which it is drawn.

Let us return to the setting on which $Y = zz^T + \sigma W$, where $W$ is a standard Wigner matrix: a symmetric matrix with iid standard gaussian entries. We want to determine for which values of $\sigma$ one excepts $X = zz^T$ to be, with high probability, the solution of \eqref{eq:thegeneralSDPformulation}, as we are interested not only to compute the MLE but also for it to coincide with the planted vector $z$ we want to recover. Since $\diag(z)W\diag(z)\sim W$ we can, without loss of generality, take $z=\1$. In that case, we are interested in understanding when 
\begin{equation}\label{eq:conditionZ2synchgaussian:11}
\lambda_2\left(D_{\left[\1\1^T + \sigma W \right]} - \left( \1\1^T + \sigma W\right) \right)>0. 
\end{equation}
Since
\[
 D_{\left[\1\1^T + \sigma W \right]} - \left( \1\1^T + \sigma W\right) = \left(nI_{n\times n} - \1\1^T\right) - \sigma\left( -D_{W} + W\right) = L_{\1\1^T} - \sigma L_{\left[-W\right]},
\]
 and $\1$ is always in the nullspace of any Laplacian matrix, it is not difficult to see that \eqref{eq:conditionZ2synchgaussian:11} is equivalent to
\begin{equation}\label{eq:conditionZ2synchgaussian:12}
 \lambda_{\max}\left( L_{\left[-W\right]} \right) < \frac{n}{\sigma}.
\end{equation}

The triangular inequality tells us that $\lambda_{\max}\left( L_{\left[-W\right]} \right) \leq \lambda_{\max}\left( -D_{W} \right) + \|W\|$. It is well known that, for any $\eps>0$, $\|W\|\leq (2+\eps)\sqrt{n}$ with high probability (see, for example, Theorem II.11 in~\cite{Szarek_Davidson_survey}). On the other hand,
\[
\lambda_{\max}\left( -D_{W} \right) = \max_{i\in [n]}\left[ -\left( D_{W} \right)_{ii} \right],
\]
which is the maximum of $n$ gaussian random variables each with variance $n$. A simple union bound yields that, for any $\eps>0$, $\lambda_{\max}\left( D_{[-W]} \right)< \sqrt{(2+\eps)n\log n}$ with high probability. This readily implies an exact recovery guarantee for \(\mathbb{Z}_2\) Synchronization with gaussian noise.

\begin{proposition}
 Let $z\in\{\pm1\}^n$ and $Y = zz^T + \sigma W$ where $W$ is a symmetric matrix with iid standard gaussian entries. If there exists $\eps>0$ such that $\sigma < \sqrt{\frac{n}{(2+\eps)\log n}}$ then, with high probability, $X = zz^T$ is the unique solution to the Semidefinite Program~\eqref{eq:thegeneralSDPformulation}.
\end{proposition}

Let us investigate the optimality of this upper bound on $\sigma$. If the diagonal elements of $D_{[-W]}$ were independent\footnote{The diagonal entries of $D_W$ are not independent because each pair of sums shares a term $W_{ij}$ as a summand.}, their distribution would be known to indeed concentrate around $\sqrt{2n \log n}$, suggesting that
\begin{equation}\label{eq:WsmallerthanDW}
 \|W\| \ll \lambda_{\max}\left( D_{[-W]} \right),
\end{equation}
which would imply
\begin{equation}\label{eq:WsmallerthanDW_L}
 \lambda_{\max}\left( L_{\left[-W\right]} \right) = \left[1+o(1)\right] \lambda_{\max}\left( D_{[-W]} \right).
\end{equation}
Both of these statements can be rigorously shown to be true. While a simple adaptation of the proof of Theorem~\ref{theorem:main} can establish~\eqref{eq:WsmallerthanDW} and~\eqref{eq:WsmallerthanDW_L} we omit their proofs for the sake of brevity, but emphasize that in this particular setting (where $W$ is a standard Wigner matrix), one does not need the whole strength of Theorem~\ref{theorem:main} as simple elementary proofs exist.

This suggests that, in rough terms, the success of the relaxation~\eqref{eq:thegeneralSDPformulation} depends mostly on whether $\lambda_{\max}\left( D_{[-W]} \right) < \frac{n}{\sigma}$, which is equivalent to
\begin{equation}\label{eq:sumWijn}
 \max_{i \in [n]}\left[ -\sigma\sum_{j=1}^{n}W_{ij} \right] < n,
\end{equation}
which can be interpreted as a bound on the amount of noise per row of $Y$. We argue next that this type of upper bound is indeed necessary for any method to succeed at recovering $z$ from $Y$.

Once again, let us consider $z=\1$ without loss of generality. Let us consider an oracle version of problem on which one is given the correct label of every single node except of node $i$. It is easy to see that the maximum likelihood estimator for $z_i$ on this oracle problem is given by
\[
 \mathrm{sign}\left[ \sum_{j \in [n]\setminus i} Y_{ij} \right] = \mathrm{sign}\left[ n-1 + \sigma\sum_{j \in [n]\setminus i} W_{ij} \right],
\]
 which would give the correct answer if and only if
\begin{equation}\label{eq:sumWijn1}
 -\sigma \sum_{j \in [n]\setminus i} W_{ij} < n -1.
\end{equation}

This means that if 
\begin{equation}\label{eq:sumWijn12}
 \max_{i \in [n]}\left[  -\sigma \sum_{j \in [n]\setminus i} W_{ij} \right] > n -1,
\end{equation}
one does not expect the MLE to succeed (with high probability) at recovering $z$ from $Y = zz^T + \sigma W$. This means that (with a uniform prior on $z$) no method is able to recover $z$ with high probability. Note the similarity between \eqref{eq:sumWijn} and \eqref{eq:sumWijn12}. This strongly suggests the optimality of the semidefinite programming based approach~\eqref{eq:thegeneralSDPformulation}.

These optimality arguments can be made rigorous. In fact, in Section~\ref{section:applications}, we will establish precise optimality results of this type, for the applications we are interested in. The main ingredient~\eqref{eq:WsmallerthanDW} in the rough argument above was the realization that the spectral norm of $W$ is, with high probability, asymptotically smaller than the largest diagonal entry of $D_{[-W]}$. Theorems~\ref{theorem:main} and~\ref{theorem:2} establish precisely this fact for a large class of matrices with independent off-diagonal entries. Empowered with this result, we will be able to establish optimality for the semidefinite programming approach to solve the problems of $\ZZ_2$ Synchronization and recovery in the stochastic block model, where the underlying random matrices have much less well understood distributions. Modulo the use of Theorem~\ref{theorem:main}, the arguments used will be very reminiscent of the the ones above.

It is pertinent to compare this approach with the one of using noncommutative Khintchine inequality, or the related matrix concentration inequalities~\cite{Tropp:TailBoundsRM,Tropp_randommatrix_monograph}, to estimate the spectral norms in question. Unfortunately, those general purpose methods are, in our case, not fine enough to give satisfactory results. One illustration of their known suboptimality is the fact that the upper bound they give for $\|W\|$ is of order $\sqrt{n\log n}$, which does not allow to establish~\eqref{eq:WsmallerthanDW}, a crucial step in the argument. In fact, the looseness of these bounds is reflected in the suboptimal guarantees obtained in~\cite{Abbe_Z2Synch,Abbe_Z2SynchER,SBM_paper_arxiv_v4}. Our results are able to establish a phenomenon of the type of~\eqref{eq:WsmallerthanDW} by relying on recent sharp estimates for the spectral norm of matrices with independent entries in~\cite{Bandeira_NARandomMatrixBound}.

\section{Main Results}\label{section:mainresults}

We use this section to formulate precise versions of, and briefly discuss, our main results.

\begin{theorem}\label{theorem:main}

Let $L$ be an $n\times n$ symmetric random Laplacian matrix (i.e. satisfying $L \1=0$) with centered independent off-diagonal entries such that
\(\sum_{j\in [n]\setminus i}\EE L_{ij}^2 \) is equal for every $i$.

Define $\sigma$ and $\sigma_{\infty}$ as 
\[
 \sigma^2 = \sum_{j\in [n]\setminus i}\EE L_{ij}^2  \quad \text{ and } \quad \sigma_{\infty}^2 = \max_{i\neq j}\left\| L_{ij} \right\|_{\infty}^2.
\]

If there exists $c>0$ such that
\begin{equation}\label{eq:req:log1}
 \sigma \geq c\left( \log n\right)^{\frac12}\sigma_{\infty},
\end{equation}
then there exists $c_1$, $C_1$, $\beta_1$, all positive and depending only on $c$, such that
\[
 \lambda_{\max}(L) \leq \left( 1 + \frac{C_1}{(\log n)^{\frac12}} \right) \max_{i}L_{ii}
\]
with probability at least $1-c_1n^{-\beta_1}$.

\end{theorem}

Even though we were not able to find a convincing application for which $\frac{\sigma}{\sigma_\infty}$ was asymptotically growing but slower than $\sqrt{\log n}$, we still include the theorem below for the sake of completeness.

\begin{theorem}\label{theorem:2}

Let $L$ be an $n\times n$ symmetric random Laplacian matrix (i.e. satisfying $L \1=0$) with centered independent off-diagonal entries such that
\(\sum_{j\in [n]\setminus i}\EE L_{ij}^2 \) is equal for every $i$.

Define $\sigma$ and $\sigma_{\infty}$ as 
\[
 \sigma^2 = \sum_{j\in [n]\setminus i}\EE L_{ij}^2  \quad \text{ and } \quad \sigma_{\infty}^2 = \max_{i\neq j}\left\| L_{ij} \right\|_{\infty}^2.
\]

If there exist $c$ and $\gamma>0$ such that
\begin{equation}\label{eq:req:log2}
 \sigma \geq c\left( \log n\right)^{\frac14 + \gamma}\sigma_{\infty},
\end{equation}
then there exist $C_2$, $c_2$, $\epsilon$ and $\beta_2$, all positive and depending only on $c$ and $\gamma>0$, such that
\[
 \lambda_{\max}(L) \leq \left( 1 + \frac{C_2}{(\log n)^{\epsilon}} \right) \max_{i}L_{ii},
\]
with probability at least $1-c_2\exp\left[ - \left( \log n\right)^{\beta_2} \right]$.

\end{theorem}

\begin{remark}
 In the theorems above, the condition that $\sum_{j\in [n]\setminus i}\EE L_{ij}^2 $ is equal for every $i$, can be relaxed to the requirement that $$ c'\sigma^2 \leq \sum_{j\in [n]\setminus i}\EE L_{ij}^2  \leq \sigma^2,$$ for all $i$. This requires only simple adaptations to the proofs of these theorems.
\end{remark}

While we defer the proof of these theorems to Section~\ref{section:proofs}, we briefly describe its idea. Lemma~\ref{corollary:Master} (borrowed from~\cite{Bandeira_NARandomMatrixBound}) estimates that
\[
 \|X\| \lesssim \sigma + \sigma_{\infty}\sqrt{\log n},
\]
where $-X$ is the off-diagonal part of $L$. One the other hand, $L_{ii} = \sum_{j\in [n]\setminus i} X_{ij}$ has variance $\sigma^2$ and the Central Limit Theorem would suggest that $L_{ii}$ behave like independent gaussians of variance $\sigma^2$, which would mean that $\max_{i}L_{ii} \sim \sigma\sqrt{\log n}$ rendering the contribution of the off-diagonal entries (to the largest eigenvalue) negligible. However, several difficulties arise: the diagonal entries are not independent (as each pair shares a summand) and one needs to make sure that the central limit theorem behavior sets in (this is, in a way, ensured by requirements \eqref{eq:req:log1} and \eqref{eq:req:log2}). The proofs in Section~\ref{section:proofs} make many needed adaptations to this argument to make it rigorous.

\section{Applications}\label{section:applications}

We now turn our attention to applications of the main results. As a form of warm-up we will start with understanding connectivity of Erd\H{o}s--R\'{e}nyi graphs.

\subsection{Connectivity of Erd\H{o}s--R\'{e}nyi graphs}\label{subsection:ErdosRenyi}

Recall that, for an integer $n$ and an edge probability parameter $0\leq p\leq 1$, the Erd\H{o}s--R\'{e}nyi graph model~\cite{ErdosRenyi_59} $\GGG(n,p)$ is a random graph on $n$ nodes where each one of the ${n \choose 2}$ edges appears independently with probability $p$. 

We are interested in understanding the probability that $G$, drawn according to $\GGG(n,p)$, is a connected graph. We will restrict our attention to the setting $p\leq \frac12$. Let $L$ be the Laplacian of the random graph, given by $D-A$ where $A$ is its adjacency matrix and $D$ a diagonal matrix containing the degree of each node. It is well-known (see, e.g.,~\cite{FanChung_SpectralGraphTheory}) that $G$ connected is equivalent to $\lambda_2(L)>0$.

It is clear that if $G$ has an isolated node then it cannot be connected. It is also known that for there not to be isolated nodes one needs the average degree of each node to be at least logarithmic~\cite{ErdosRenyi_59}. For this reason we will focus on the regime
\[
 p = \frac{\rho \log n}{n},
\]
for a constant $\rho$. It is easy to establish a phase transition on the degrees of the nodes of graphs drawn from $\GGG(n,p)$.

\begin{lemma}\label{ER_phasetransition_degrees}

Let $n$ be a positive integer, $\rho$ a constant, and $p =  \frac{\rho \log n}{n}$. Let $G$ be a random graph drawn from  $\GGG(n,p)$, then for any constant $\Delta>0$:
\begin{enumerate}
 \item If $\rho > 1$ then, with high probability, $\min_{i\in [n]}\deg(i) \geq\frac{ \Delta }{\sqrt{\log n}}\EE \deg(i)$.
 \item If $\rho < 1$ then, with high probability, $\min_{i\in [n]}\deg(i) = 0$. That is, $G$ has at least one isolated node, thus being disconnected.
\end{enumerate}

\end{lemma}

Part (2) of the Lemma is a classical result~\cite{ErdosRenyi_59}, a particularly simple proof of it proceeds by applying the second moment method to the number of isolated nodes in $G$. For the sake of brevity we will skip those details, and focus on part (1). The main thing to note in part (1) of Lemma~\ref{ER_phasetransition_degrees} is that the lower bound on minimum degree is asymptotically smaller than the average degree $\EE\deg(i)$.

\proof{[of part (1) of Lemma~\ref{ER_phasetransition_degrees}]

Let $p=\frac{\rho \log n}{n}$ and $i$ denote a node of the graph, note that $\EE\deg(i) = \frac{n-1}n\rho \log n$. We use Chernoff bound (see, for example, Lemma 2.3.3 in~\cite{Durrett_RandomGraphDynamics}) to establish, for any $0<t<1$,
\begin{eqnarray*}
  \PP\left[ \deg(i) < t \EE\deg(i)  \right] & \leq & \left[ \frac{\exp(-(1-t))}{t^t} \right]^{\EE\deg(i)} \\
& = & \left[ \frac{\exp(-(1-t))}{t^t} \right]^{\frac{n-1}{n}\rho\log n} \\
& = & \exp\left[-\left[1-t-t\log(1/t)\right]\frac{n-1}{n}\rho\log n\right].
\end{eqnarray*}

Taking $t = \frac{\Delta}{\sqrt{\log n}}$ gives, for $n$ large enough (so that $t\leq 1$), that the probability that $\deg(i) < \frac{\Delta}{\sqrt{\log n}} \EE\deg(i)$ is at most 
\[
 \exp\left[-\left[1-\frac{\Delta}{\sqrt{\log n}}-\frac{\Delta}{\sqrt{\log n}}\log\left(\frac{\sqrt{\log n}}{\Delta}\right)\right]\frac{n-1}{n}\rho\log n\right],
\]
which is easily seen to be $\exp\left[ -\rho\log n + O(\sqrt{\log n} \log\log n)   \right]$. A simple union bound over the $n$ vertices of $G$ gives
\[
 \PP\left[ \min_{i\in [n]}\deg(i) <  \frac{\Delta}{\sqrt{\log n}} \EE\deg(i)  \right] \leq  \exp\left[ -(\rho-1)\log n + O(\sqrt{\log n} \log\log n)   \right].
\]
\qedhere

}

Using Theorem~\ref{theorem:main} we will show that, with high probability, as long as every node in $G$ is at least $\frac{\Delta}{\sqrt{\log n}}$ of the average degree, for a suitable $\Delta$, then $G$ is connected. This is made precise in the following Lemma.

\begin{lemma}\label{ER_DegreeVsConnectivity}
 Let $n\geq 2$ be an integer and $\eps>0$. Suppose that $ \frac{\eps \log n}{n}\leq p \leq \frac12$ and $G$ a random graph drawn from $\GGG(n,p)$. There exists a constant $\Delta$ such that, with high probability, the following holds:

If
\[
 \min_{i\in [n]}\deg(i) \geq \frac{\Delta}{\sqrt{\log n}} \EE \deg(i),
\]
then $G$ is a connected graph (note that the right hand side does not depend on $i$).
\end{lemma}

Before proving this Lemma, we note that Lemmas~\ref{ER_phasetransition_degrees} and~\ref{ER_DegreeVsConnectivity} immediately imply the well known phase transition phenomenon.

\begin{theorem}\label{ER_phasetransition}
Let $n$ be a positive integer and $p =  \frac{\rho \log n}{n}$.
\begin{enumerate}
 \item If $\rho > 1$ then, with high probability, a random graph drawn from $\GGG(n,p)$ is connected.
 \item If $\rho < 1$ then, with high probability, a random graph drawn from $\GGG(n,p)$ has at least one isolated node, thus being disconnected.
\end{enumerate}

\end{theorem}

While this phase transition is well understood, we find our proof through Lemmas~\ref{ER_phasetransition_degrees} and~\ref{ER_DegreeVsConnectivity} enlightening, as it provides a simple explanation of why the phase transition for disappearance of isolated nodes coincides with the one for connectivity. Moreover, it  also emphasizes a connection with the optimality of the semidefinite relaxations in both $\ZZ_2$ Synchronization and the Stochastic Block Model that we will discuss in the sections to follow.

\proof{[of Lemma~\ref{ER_DegreeVsConnectivity}]

Let $L$ be the graph Laplacian of $G$. Note that $\EE(L) = npI - p\1\1^T$, which means that
\[
 L = npI - p\1\1^T - \left[-L+\EE(L)\right]
\]

Since $L\1 = 0$, it is easy to see that $G$ is connected if and only if

\[
 \lambda_{\max}\left[-L+\EE(L)\right] < np
\]

We proceed by using Theorem~\ref{theorem:main} for

\[
 L = -L+\EE(L).
\]

The hypotheses of the Theorem are satisfied as the off-diagonal entries of $L$ are independent and
\[
   \sum_{j\in [n]\setminus i}\EE L_{ij}^2 = (n-1)p(1-p) \geq \frac{np(1-p)}2 \geq \frac{\eps}2 (1-p)^2 \log n = \frac{\eps}2 \log n \max_{i\neq j}\left\| L_{ij} \right\|_\infty^2.
\]
This guarantees that there exists a constant $C_1$ such that, with high probability, 

\begin{equation}\label{eq:eventofconnectivityform}
\lambda_{\max}\left[-L+\EE(L)\right] \leq \left(1+\frac{C_1}{\sqrt{\log n}}\right)\max_{i\in [n]}\left[ - \deg(i) + (n-1)p  \right]
\end{equation}
where $\deg(i) = L_{ii}$ is the degree of node $i$. Equivalently,
\begin{eqnarray*}
\lambda_{\max}\left[-L+\EE(L)\right] & \leq & np + \left(1+\frac{C_1}{\sqrt{\log n}}\right)\left[ - \min_{i\in [n]}\deg(i) + (n-1)p  \right] - np
\end{eqnarray*}

This means that, as long as \eqref{eq:eventofconnectivityform} holds, then
\[
 \left(1+\frac{C_1}{\sqrt{\log n}}\right)\left[ - \min_{i\in [n]}\deg(i) + (n-1)p  \right] - np< 0 
\]
implies the connectivity of $G$. Straighforward manipulations show that this condition is equivalent to

%
%
%
%
%
%
%
%
%

\begin{equation*}
 \min_{i} \deg(i) > np\frac{C_1}{\sqrt{\log n}+C_1} - p,
\end{equation*}

which is implied by

\begin{equation}
 \min_{i} \deg(i) \geq np\frac{C_1}{\sqrt{\log n}}.
\end{equation}

The lemma follows by taking $\Delta = 2C_1$.

\qedhere

}

%
%
%
%
%

\subsection{Synchronization over the group of two elements}\label{subsection:Int:Z2Synch}

Recall the setting of $\ZZ_2$ Synchronization~\cite{Abbe_Z2Synch,Abbe_Z2SynchER}. Given an underlying graph $G$ with $n$ nodes, the task is to recover a binary vector $z\in\{\pm 1\}^n$ from noisy measurements $Y_{ij}$ of $z_iz_j$. Following~\cite{Abbe_Z2Synch,Abbe_Z2SynchER} we will take the underlying graph $G$ to be an Erd\H{o}s--R\'{e}nyi graph $\GGG(n,p)$ and, for each edge $(i,j)\in G$,
\[
Y_{ij} = 
\left\{
\begin{array}{rcl}
z_iz_j & \text{ with probability } & 1-\eps \\
-z_iz_j & \text{ with probability } & \eps,
\end{array}
\right.
\]
 where $\eps<\frac12$ represents the noise level. We are interested in understanding for which values of $p$ and $\eps$ is it possible to exactly recover $z$. It is easy to see that, just like in the example in Section~\ref{section:firstexample}, the maximum likelihood estimator is given by~\eqref{littleGrothendieck:Z2:1}. Similarly, we consider its semidefinite relaxation \eqref{eq:thegeneralSDPformulation} and investigate when $X = zz^T$ is the unique solution of \eqref{eq:thegeneralSDPformulation}.

It is easy to see that $Y$ is given by
\[
 Y = \diag(z)\left(A_G - 2A_H\right)\diag(z),
\]
where $A_G$ is the adjacency matrix of the underlying graph and $A_H$ is the adjacency of the graph consisting of the corrupted edges. In this case we want conditions on $\eps$ and $p$ under which $zz^T$ is the unique solution to:

\begin{equation}\label{SDP:forZ2Synch}
\begin{array}{cl}
 \max & \Tr\left[ \diag(z)\left( A_G - 2A_H \right)\diag(z)X \right] \\
 \text{ s.t. }  & X_{ii} = 1\\
		& X\succeq 0.
\end{array}
\end{equation}

Lemma~\ref{lemma:dualcertificate_constructed} states that $zz^T$ is indeed the unique solution as long as the second smallest eigenvalue of
\begin{equation}\label{eq:conditioneig:143}
D_{A_G - 2A_H} -\diag(z)\left(A_G-2A_H\right)\diag(z)= D_G - 2D_H - \diag(z)\left( A_G - 2A_H \right)\diag(z) 
\end{equation}
is strictly positive. As $\diag(z)\left( D_G - 2D_H \right)\diag(z) = D_G - 2D_H$ and conjugating by $\diag(z)$ does not alter the eigenvalues, the second smallest eigenvalue of \eqref{eq:conditioneig:143} being strictly positive is equivalent to
\begin{equation}\label{eq:conditioneig:144}
 \lambda_2\left( D_G - A_G - 2\left( D_H - A_H \right)   \right) >0.
\end{equation}

Since $D_G - A_G - 2\left( D_H - A_H \right) = L_G - 2L_H$, where $L_G$ and $L_H$ are the Laplacians of, respectively, $G$ and $H$, we define $L_{\mathrm{Synch}}$ and write the condition in terms of $L_{\mathrm{Synch}}$.

\begin{definition}\label{def:Lsynch}\emph{[$L_{\mathrm{Synch}}$]}
In the setting described above,
 \[
  L_{\mathrm{Synch}} = L_G - 2L_H,
 \]
 where $G$ is the graph of all measurements and $H$ is the graph of wrong measurements.
\end{definition}

Then, \eqref{eq:conditioneig:144} is equivalent to $\lambda_2\left(L_{\mathrm{Synch}}\right) >0$. The following Lemma readily follows by noting that $\EE \left[ L_{\mathrm{Synch}} \right] = np(1-2\eps)I_{n\times n} - p(1-2\eps)\1\1^T$.

\begin{lemma}\label{lemma:conditionSDPistightZ2Synch}
Consider the $\ZZ_2$ Synchronization problem defined above and $L_{\mathrm{Synch}}$ defined in Definition~\ref{def:Lsynch}. As long as 
\[
 \lambda_{\max} \left( - L_{\mathrm{Synch}} + \EE \left[ L_{\mathrm{Synch}} \right] \right) < np(1-2\eps),
\]
the Semidefinite program~\eqref{SDP:forZ2Synch} achieves exact recovery.
\end{lemma}

In~\cite{Abbe_Z2Synch,Abbe_Z2SynchER}, this largest eigenvalue is estimated using the general purpose matrix concentration inequalities (such as the ones in~\cite{Tropp:TailBoundsRM}) obtaining a suboptimal bound. In contrast, we will do this estimate using Theorem~\ref{theorem:main}.

Let us define, for a node $i$, $\deg_+(i)$ as the number of non-corrupted edges incident to $i$ and $\deg_-(i)$ as the number of corrupted edges incident to $i$.
We start by obtaining the following theorem.

\begin{theorem}\label{Theorem:Z2SynchDelta}
As long as $n>2$, $p>\frac{\log n}{2n}$ and $p(1-2\eps)^2 \leq \frac12$, there exists $\Delta>0$ such that, with high probability, the following holds: If
\begin{equation}\label{eq:inTheorem:Z2SynchDelta}
\min_{i\in [n]}\left[\deg_{+}(i) - \deg_{-}(i)\right] \geq \frac{\Delta}{\sqrt{\log n}} \EE \left[ \deg_{+}(i) - \deg_{-}(i) \right],
\end{equation}
then the semidefinite program \eqref{SDP:forZ2Synch} achieves exact recovery.
\end{theorem}


\proof{[of Theorem~\ref{Theorem:Z2SynchDelta}]

The idea is to apply Theorem~\ref{theorem:main} to $L = - L_{\mathrm{Synch}} + \EE \left[ L_{\mathrm{Synch}}  \right] $. Note that $L$ has independent off-diagonal entries and
\begin{eqnarray*}
\sum_{j\in [n]\setminus i}\EE\left[L_{ij}^2\right] & = & (n-1)\left( p-p^2(1-2\eps)^2 \right)
\geq \frac14 np \geq \frac18 \log n 
\\ & \geq & \frac{1+ p(1-2\eps)}{8(1+\sqrt{2})} \log n
=  \frac{\log n}{8(1+\sqrt{2})} \max_{i\neq j}\left\|L_{ij}^2\right\|_\infty.
\end{eqnarray*}

Hence, there exists a constant $\Delta'$ such that, with high probability,
\[
 \lambda_{max}\left(  - L_{\mathrm{Synch}} + \EE \left[ L_{\mathrm{Synch}} \right] \right) \leq \left(1+\frac{\Delta'}{\sqrt{\log n}} \right) \max_{i\in [n]} \left[  - (L_{\mathrm{Synch}})_{ii} + \EE \left[ (L_{\mathrm{Synch}})_{ii} \right] \right].
\]

We just need to show that, there exists $\Delta>0$ such that, if \eqref{eq:inTheorem:Z2SynchDelta} holds, then
\begin{equation}\label{eq:whatneedstobeprovedZ2Synch32}
 \left(1+\frac{\Delta'}{\sqrt{\log n}} \right) \max_{i\in [n]} \left[  - (L_{\mathrm{Synch}})_{ii} + \EE \left[ (L_{\mathrm{Synch}})_{ii} \right] \right] < np(1-2\eps).
\end{equation}

Recall that $ (L_{\mathrm{Synch}})_{ii} = \deg_{+}(i) - \deg_{-}(i) $ and $\EE (L_{\mathrm{Synch}})_{ii} = (n-1)p(1-2\eps)$. We can rewrite \eqref{eq:whatneedstobeprovedZ2Synch32} as
\[
  \min_{i\in [n]}(L_{\mathrm{Synch}})_{ii}     > (n-1)p(1-2\eps) -  np(1-2\eps)\left(1+\frac{\Delta'}{\sqrt{\log n}} \right)^{-1}.
\]

Straightforward algebraic manipulations show that there exists a constant $\Delta$ such that
\[
 (n-1)p(1-2\eps) -  np(1-2\eps)\left(1+\frac{\Delta'}{\sqrt{\log n}} \right)^{-1} \leq \frac{\Delta}{\sqrt{\log n}} \EE \left[ \deg_{+}(i) - \deg_{-}(i) \right],
\]
proving the Theorem.

\qedhere

}


We note that, if $p \leq \frac{\log n}{2n}$, then Theorem~\ref{ER_phasetransition} implies that, with high probability, the underlying graph is disconnected implying impossibility of exact recovery. We also note that if we do not have
\begin{equation}\label{eq:inTheorem:Z2SynchDelta0}
 \min_{i\in [n]}\left[\deg_{+}(i) - \deg_{-}(i)\right] \geq 0,
\end{equation}
then the maximum likelihood does not match the ground truth, rendering exact recovery unrealistic\footnote{Recall that, if we assume a uniform prior, the MLE is the method that maximizes the probability of exact recovery}. The optimality of this analysis hinges upon the fact that the right-hand side of \eqref{eq:inTheorem:Z2SynchDelta} is asymptotically smaller than the expectation of $\deg_{+}(i) - \deg_{-}(i)$, suggesting that \eqref{eq:inTheorem:Z2SynchDelta} and \eqref{eq:inTheorem:Z2SynchDelta0} have similar probabilities and the same phase transition.

The next Theorem establishes the optimality of the semidefinite programming based approach in a particular regime, solving a problem raised in~\cite{Abbe_Z2Synch,Abbe_Z2SynchER}. While it is clear that one can use Theorem~\ref{Theorem:Z2SynchDelta} to establish similar results for many other regimes (for some, through estimates similar to the ones in Lemma~\ref{lemma:controlingTforSBM}), the main purpose of this paper is not to perform a detailed analysis of this problem but rather to illustrate the efficacy of these semidefinite relaxations and the fundamental connections between these different phenomena, through Theorem~\ref{theorem:main}. The independent parallel research efforts of Hajek et al.~\cite{Hajek_et_al_SBM_SDP_extensions} address other regimes for this particular problem, we refer the interested reader there.

\begin{corollary}\label{corollary:Z2synchSDPisoptimal}

As long as $\eps < \frac12$ and $p(1-2\eps)^2 \leq \frac12$, there exists a constant $K$ for which the following holds: If there exists $\delta>0$ such that 
\begin{equation}\label{conditionfor:corollary:Z2synchSDPisoptimal}
(n-1)p \geq (1+\delta) \frac2{(1-2\eps)^2} \left[ 1 + \frac{K}{\sqrt{\log n}} + \frac53(1-2\eps)  \right] \log n,
\end{equation}
then the Semidefinite program~\eqref{SDP:forZ2Synch} achieves exact recovery with high probability.

\end{corollary}

Before proving this corollary we emphasize how it solves the problem, raised in~\cite{Abbe_Z2Synch,Abbe_Z2SynchER}, of whether the semidefinite programming approach for $\ZZ_2$ Synchronization is optimal in the low signal-to-noise regime. In fact, the results in~\cite{Abbe_Z2Synch,Abbe_Z2SynchER} ensure that the threshold in Corollary~\ref{corollary:Z2synchSDPisoptimal} is optimal for, at least, an interesting range of values of $\eps$. Empowered with Theorem~\ref{Theorem:Z2SynchDelta}, the proof of this corollary becomes rather elementary.

\proof{[of Corollary~\ref{corollary:Z2synchSDPisoptimal}]

This corollary will be established with a simple use of Bernstein's inequality. 

Our goal is to show that, given $\Delta$, there exists a $K$ and $\delta$ such that, under the hypothesis of the Corollary,
\[
\min_{i\in [n]}\left[\deg_{+}(i) - \deg_{-}(i)\right] \geq \frac{\Delta}{\sqrt{\log n}} \EE \left[ \deg_{+}(i) - \deg_{-}(i) \right],
\]
holds with high probability. This implies, via Theorem~\ref{Theorem:Z2SynchDelta}, that the semidefinite program~\eqref{SDP:forZ2Synch} achieves exact recovery with high probability.

We will consider $n$ to be large enough. We start by noting that it suffices to show that there exists $ \delta > 0 $ such that, for each $i\in [n]$ separately,
\begin{equation}\label{eq:whatweshowforeachiZ2}
\PP\left[   \deg_{+}(i) - \deg_{-}(i) <   \frac{\Delta}{\sqrt{\log n}} \EE \left[ \deg_{+}(i) - \deg_{-}(i) \right]  \right] \leq n^{-(1+\delta)}.
\end{equation}
Indeed, \eqref{eq:whatweshowforeachiZ2} together with a union bound over the $n$ nodes of the graph would establish the Corollary.

Throughout the rest of the proof we will fix $i\in [n]$ and use  $\deg_{+}$ and $\deg_{-}$ to denote, respectively, $\deg_{+}(i)$ and $\deg_{-}(i)$. It is easy to see that
\[
 \deg_{+} - \deg_{-} = (n-1)p(1-2\eps) -  \sum_{j=1}^{n-1} x_j,
\]
where $x_j$ are i.i.d. centered random variables with distribution
\[
x_j = 
\left\{
\begin{array}{rcl}
-1 + p(1-2\eps) & \text{ with probability } & p (1-\eps) \\
1 + p(1-2\eps) & \text{ with probability } & p \eps \\
 p(1-2\eps)  & \text{ with probability } & 1-p. \\
\end{array}
\right.
\]

For any $t>0$ Bernstein's inequality gives
\[
 \PP\left[ \sum_{j=1}^{n-1}x_j > t \right] \leq \exp\left( -\frac{t^2/2}{(n-1)\EE x_j^2 + \frac{t}3\|x_j\|_{\infty}} \right).
\]

Taking $t = \left[ 1 - \frac{\Delta}{\sqrt{\log n}}  \right] (n-1) p (1-2\eps) $ gives
\begin{eqnarray*}
 \PP\left[ \deg_{+} - \deg_{-} <   \frac{\Delta}{\sqrt{\log n}} \EE \left[ \deg_{+} - \deg_{-} \right]  \right] \\
\leq \exp\left( -\frac{\left( \left[ 1 - \frac{\Delta}{\sqrt{\log n}}  \right] (n-1) p (1-2\eps) \right)^2/2}{(n-1)\EE x_j^2 + \frac{\left( \left[ 1 - \frac{\Delta}{\sqrt{\log n}}  \right] (n-1) p (1-2\eps) \right)}3\|x_j\|_{\infty}} \right) \\
= \exp\left( -\frac{ \left[ 1 - \frac{\Delta}{\sqrt{\log n}}  \right]^2 (n-1) p (1-2\eps)^2/2}{ \frac1p\EE x_j^2 + \frac{\left( \left[ 1 - \frac{\Delta}{\sqrt{\log n}}  \right] (1-2\eps) \right)}3\|x_j\|_{\infty}} \right)
\end{eqnarray*}

Condition \eqref{conditionfor:corollary:Z2synchSDPisoptimal} (for a $K$ to be determined later) guarantees that
\[
 (n-1) p (1-2\eps)^2/2 \geq  (1+\delta) \left[ 1 + \frac{K}{\sqrt{\log n}} + \frac53(1-2\eps)  \right] \log n,
\]
meaning that we just need to show that there exists $K>0$ for which
\[
\frac{ \left[ 1 - \frac{\Delta}{\sqrt{\log n}}  \right]^2  \left(1 + \frac{K}{\sqrt{\log n}} + \frac53(1-2\eps)\right)}{ \frac1p\EE x_j^2 + \frac{\left( \left[ 1 - \frac{\Delta}{\sqrt{\log n}}  \right] (1-2\eps) \right)}3\|x_j\|_{\infty}}
  \geq 1.
\]
Note that $\frac1p\EE x_j^2 = 1 + p(1-2\eps) \leq 1 + (1-2\eps)$ and $\|x_j\|_{\infty} = 1 + p(1-2\eps) \leq 2$, implying that
\[
 \frac1p\EE x_j^2 + \frac{\left( \left[ 1 - \frac{\Delta}{\sqrt{\log n}}  \right] (1-2\eps) \right)}3\|x_j\|_{\infty}
\leq 1 + \frac53(1-2\eps).
\]
 Also, $\left[ 1 - \frac{\Delta}{\sqrt{\log n}}  \right]^{2} \geq 1 - \frac{2\Delta}{\sqrt{\log n}}$. The corollary is then proved by noting that there exists $K>0$ such that
\[
  \frac{K}{\sqrt{\log n}} \geq  2 K\frac{\Delta}{\log n} + \frac{2\Delta}{\sqrt{\log n}}\left(    1 + \frac53(1-2\eps) \right).
\]

\qedhere

}

\subsection{Stochastic Block Model with two communities}\label{subsection:Int:SBM}

We shift our attention to the problem of exact recovery of the stochastic block model with two communities. Recall Definition~\ref{def:SBM}, for $n$ even and $0\leq q<p \leq 1$, we say that a graph $G$ with $n$ nodes is drawn from the Stochastic block model with two communities $\GGG(n,p,q)$ if the nodes are divided in two sets of $\frac{n}2$ nodes each, and for each pair of vertices $i,j$, $(i,j)$ is an edge of $G$ with probability $p$ if $i$ and $j$ are in the same cluster and $q$ otherwise, independently from any other edge. Let $g\in \{\pm 1\}^{n}$ be a vector that is $1$ in one of the clusters and $-1$ in the other, our task is to recover $g$.

The maximum likelihood estimator for $g$ is given by
\begin{equation}\label{littleGrothendieck:Z2:1SBM}
\begin{array}{cl}
 \max & x^TBx \\
 \text{ s.t. } & x\in\RR^n \\
		& x_i^2 = 1, \\
		& \sum_{i=1}^n x_i = 0,
\end{array}
\end{equation}
where $B$ is the signed adjacency of $G$, meaning that $B_{ij}=1$ if $(i,j)$ is an edge of $G$ and $B_{ij}=-1$ otherwise. Note that $B = 2A - \left(\1\1^T - I \right)$, where $A$ is the adjacency matrix. We will drop the balanced constraint $\sum_{i=1}^n x_i = 0$, arriving at~\eqref{littleGrothendieck:Z2:1} for $Y=B$. The intuitive justification is that there are enough $-1$ entries in $B$ to discourage unbalanced solutions. As in the problems considered above, we will consider the semidefinite relaxation~\eqref{eq:thegeneralSDPformulation}.
\begin{equation}\label{SDP:forSBM}
\begin{array}{cl}
 \max & \Tr\left[\left( 2A - \left(\1\1^T - I \right) \right)X\right] \\
 \text{ s.t. }  & X_{ii} = 1\\
		& X\succeq 0.
\end{array}
\end{equation}

We want to understand when is it that $X=gg^T$ is the unique solution of~\eqref{SDP:forSBM}.
Lemma~\ref{lemma:dualcertificate_constructed} shows that $gg^T$ is indeed the unique solution of~\eqref{SDP:forSBM} as long as the second smallest eigenvalue of
\begin{equation}\label{eq:conditioneig:1435}
D_{\left[\diag(g)( 2A - \left(\1\1^T - I \right) )\diag(g)\right]} - \left[2A - \left(\1\1^T - I \right)\right],
\end{equation}
is strictly positive.

Let us introduce a new matrix.
\begin{definition}\emph{[$\Gamma_{\mathrm{SBM}}$]}
Given a graph $G$ drawn from the stochastic block model with two clusters, 
 \[
  \Gamma_{\mathrm{SBM}} = \DDD_+ - \DDD_- - A,
 \]
where $\DDD_+$ is a diagonal matrix of inner degrees, $\DDD_-$ is a diagonal matrix of outer degrees and $A$ is the adjacency matrix of the graph.
\end{definition}

It is easy to see that $D_{\left[\diag(g)A\diag(g)\right]} =  \DDD_+ - \DDD_-$. In fact, 
\[
D_{\left[\diag(g)( 2A - \left(\1\1^T - I \right) )\diag(g)\right]} - \left[2A - \left(\1\1^T - I \right)\right] = 2 \Gamma_{\mathrm{SBM}} + \1\1^T,
\]
which means that $gg^T$ is the unique solution of~\eqref{SDP:forSBM} as long as
$
\lambda_2\left( 2\Gamma_{\mathrm{SBM}} + \1\1^T \right) > 0.
$

Note that
\begin{eqnarray*}
\EE \left[ 2\Gamma_{\mathrm{SBM}} + \1\1^T \right] & = & 2\left(\left( \frac{n}2 p - \frac{n}2q \right)I_{n\times n} - \left( \frac{p+q}2\1\1^T + \frac{p-q}2gg^T  \right) \right) + \1\1^T \\ 
& = & n \left( p-q \right)\left(I_{n\times n} - \frac{gg^T}n \right) + n\left( 1 - (p+q) \right)\frac{\1\1^T}n.
\end{eqnarray*}

If we suppose that $p<\frac12$, we have $1 - (p+q) > p-q$ the second smallest eigenvalue of $\EE \left[ 2\Gamma_{\mathrm{SBM}} + \1\1^T \right] $ is $n \left( p-q \right)$. This establishes the following Lemma.

\begin{lemma}\label{lemma:conditionSDPistightSBM}
Let $n\geq 4$ be even and let $G$ be drawn from $G(n,p,q)$ with edge probabilities $p<\frac12$ and $q<p$. As long as 
\[
 \lambda_{\max} \left( - \Gamma_{\mathrm{SBM}} + \EE \left[ \Gamma_{\mathrm{SBM}} \right] \right) < \frac{n}2(p-q),
\]
the Semidefinite program~\eqref{SDP:forSBM} for the stochastic block model problem achieves exact recovery, meaning that $gg^T$ is its unique solution.
\end{lemma}

Estimating this largest eigenvalue using Theorem~\ref{theorem:main}, we obtain the following theorem.

\begin{theorem}\label{Theorem:SBMDelta}
Let $n\geq 4$ be even and let $G$ be drawn from $G(n,p,q)$. As long as $\frac{\log n}{3n}<p<\frac12$ and $q<p$, then there exists $\Delta>0$ such that, with high probability, the following holds: If,
\begin{equation}\label{eq:inTheorem:SBMDelta}
\min_{i}\left( \deg_{in}(i) - \deg_{out}(i) \right) \geq \frac{\Delta}{\sqrt{\log n}} \EE \left[ \deg_{in}(i) - \deg_{out}(i) \right]
\end{equation}
then the semidefinite program \eqref{SDP:forSBM} achieves exact recovery.
\end{theorem}

\proof{

The idea is again to apply Theorem~\ref{theorem:main}. One obstacle is that $\Gamma_{\mathrm{SBM}}$ is not a Laplacian matrix. Let $g$ denote the vector that is $1$ in a cluster and $-1$ in the other, and let $\diag(g)$ denote a diagonal matrix with the entries of $g$ on the diagonal. We define
\[
 \Gamma'_{\mathrm{SBM}} = \diag(g)  \Gamma_{\mathrm{SBM}} \diag(g).
\]
Note that $\Gamma'_{\mathrm{SBM}}$ is a Laplacian and both the eigenvalues and diagonal elements of $  \EE \left[ \Gamma'_{\mathrm{SBM}} \right] - \Gamma'_{\mathrm{SBM}}$ are the as $  \EE \left[ \Gamma_{\mathrm{SBM}} \right] - \Gamma_{\mathrm{SBM}} $.

We apply Theorem~\ref{theorem:main} to $L = - \Gamma'_{\mathrm{SBM}} + \EE \left[ \Gamma'_{\mathrm{SBM}} \right] $. Note that $L$ has independent off-diagonal entries and
\begin{eqnarray*}
\sum_{j\in [n]\setminus i}\EE\left[L_{ij}^2\right] & = & \left(\frac{n}2-1\right)\left(p-p^2\right) + \frac{n}2\left(q-q^2\right)
\geq \frac{n}8p \geq \frac{\log n}{24}
\\ & \geq & 
 \frac{\log n}{24}(1-q)  =   \frac{\log n}{24}\max_{i\neq j}\left\|L_{ij}^2\right\|_\infty.
\end{eqnarray*}

Hence, there exists a constant $\Delta'$ such that, with high probability,
\[
 \lambda_{max}\left(   - \Gamma'_{\mathrm{SBM}} + \EE \left[ \Gamma'_{\mathrm{SBM}} \right] \right) \leq \left(1+\frac{\Delta'}{\sqrt{\log n}} \right) \max_{i\in [n]} \left[  - (\Gamma'_{\mathrm{SBM}})_{ii} + \EE \left[ (\Gamma'_{\mathrm{SBM}})_{ii} \right] \right],
\]

which is equivalent to
\begin{equation}
 \lambda_{max}\left(   - \Gamma_{\mathrm{SBM}} + \EE \left[ \Gamma_{\mathrm{SBM}} \right] \right) \leq \left(1+\frac{\Delta'}{\sqrt{\log n}} \right) \max_{i\in [n]} \left[  - (\Gamma_{\mathrm{SBM}})_{ii} + \EE \left[ (\Gamma_{\mathrm{SBM}})_{ii} \right] \right].
\end{equation}

We just need to show that, there exists $\Delta>0$ such that, if \eqref{eq:inTheorem:SBMDelta} holds, then
\begin{equation}\label{eq:whatneedstobeprovedSBM32}
 \left(1+\frac{\Delta'}{\sqrt{\log n}} \right) \max_{i\in [n]} \left[  - (\Gamma_{\mathrm{SBM}})_{ii} + \EE \left[ (\Gamma_{\mathrm{SBM}})_{ii} \right] \right] < \frac{n}2(p-q)-p.
\end{equation}

Note that $(\Gamma_{\mathrm{SBM}})_{ii} = \deg_{in}(i) - \deg_{out}(i)$ and 
\[
\EE \left[ \deg_{in}(i) - \deg_{out}(i) \right] = \frac{n}2(p-q) - p. 
\]

Condition \eqref{eq:inTheorem:SBMDelta} can thus be rewriten as
\[
 \max_{i\in [n]} \left[  - (\Gamma_{\mathrm{SBM}})_{ii} + \EE \left[ (\Gamma_{\mathrm{SBM}})_{ii} \right] \right]  \leq \left[ 1- \frac{\Delta}{\sqrt{\log n}}\right] \left( \frac{n}2(p-q) - p \right) .
\]

The Theorem is then proven by noting that, for any $\Delta'$, there exists $\Delta$ such that
\[
 \left[ 1- \frac{\Delta}{\sqrt{\log n}}\right] \left( \frac{n}2(p-q) - p \right) \leq \left[1+\frac{\Delta'}{\sqrt{\log n}} \right]^{-1} \left( \frac{n}2(p-q)-p\right).
\]

\qedhere

}

As a corollary of this theorem we can establish a sharp threshold for exact recovery for the stochastic block model of two clusters solving a problem posed in~\cite{SBM_paper_arxiv_v4}. We recall that this problem was simultaneously solved by the parallel research efforts of Hajek et al.~\cite{Hajek_et_al_SBM_SDP}. 

We first show a Lemma concerning $\min_{i}\left( \deg_{in}(i) - \deg_{out}(i) \right)$, analogous to Lemma~\ref{ER_phasetransition_degrees}.

\begin{lemma}\label{corollary:SBMthreshold:onlyaboutdegrees}
Let $G$ be a random graph with $n$ nodes drawn accordingly to the stochastic block model on two communities with edge probabilities $p$ and $q$. Let $p = \frac{\alpha \log n}n$ and $q = \frac{\beta \log n}n$, where $\alpha> \beta$ are constants. Then for any constant $\Delta > 0$,
\begin{enumerate}
\item If
\begin{equation}\label{eq:conditionsalphabetaSBMcor}
 \sqrt{\alpha} - \sqrt{\beta} > \sqrt{2},
\end{equation}
then, with high probability,
\[
\min_{i}\left( \deg_{in}(i) - \deg_{out}(i) \right) \geq \frac{\Delta}{\sqrt{\log n}} \EE \left[ \deg_{in}(i) - \deg_{out}(i) \right].
\]

\item On the other hand, if 
\begin{equation}\label{eq:conditionsalphabetaSBMcor_no}
 \sqrt{\alpha} - \sqrt{\beta} < \sqrt{2},
\end{equation}
then, with high probability,
\[
\min_{i}\left( \deg_{in}(i) - \deg_{out}(i) \right) < 0,
\]
and exact recovery is impossible.
\end{enumerate}
\end{lemma}

Part (2) is proven in~\cite{SBM_paper_arxiv_v4}, so we will focus on part (1). Before proving this lemma we note how, together with Theorem~\ref{Theorem:SBMDelta}, this immediately implies the following Corollary.

\begin{corollary}\label{corollary:SBMthreshold}

Let $G$ be a random graph with $n$ nodes drawn accordingly to the stochastic block model on two communities with edge probabilities $p$ and $q$. Let $p = \frac{\alpha \log n}n$ and $q = \frac{\beta \log n}n$, where $\alpha> \beta$ are constants. Then, as long as
\begin{equation}\label{eq:conditionsalphabetaSBMcor_2}
 \sqrt{\alpha} - \sqrt{\beta} > \sqrt{2},
\end{equation}
 the semidefinite program~\eqref{SDP:forSBM} coincides with the true partition with high probability.  

\end{corollary}


In order to establish Lemma~\ref{corollary:SBMthreshold:onlyaboutdegrees} we will borrow an estimate from~\cite{SBM_paper_arxiv_v4}.
 
\begin{definition}\label{def:definitionofT}\emph{[Definition 3 in~\cite{SBM_paper_arxiv_v4}]}
Let $m$ be a natural number, $p,q\in [0,1]$, and $\delta\in \RR$, we define
\[
T(m, p, q , \delta) = \PP\left[ \sum_{i=1}^m (Z_i - W_i) \geq \delta \right],
\]
where $W_1,\dots, W_m$ are i.i.d.\ $\mathrm{Bernoulli}(p)$ and $Z_1,\dots, Z_m$ are i.i.d.\ $\mathrm{Bernoulli}(q)$, independent of $W_1,\dots, W_m$.
\end{definition}

\begin{lemma}\label{lemma:controlingTforSBM}
Recall Definition~\ref{def:definitionofT}. Let $\alpha$, $\beta$, and $\Delta'$ be constants. Then,
 \[
 T\left(\frac{n}2,\frac{\alpha\log n}{n},\frac{\beta\log n}{n},-\Delta' \sqrt{\log n}\right) \leq \exp\left[ -\left(\frac{\alpha + \beta}2 - \sqrt{\alpha\beta} - \delta(n)\right)\log n\right],
 \]
with $\displaystyle{\lim_{n \to \infty}\delta(n) = 0}$.
\end{lemma}

\proof{  

The proof of this Lemma is obtained by straightforward adaptations to the proof of Lemma 8 in~\cite{SBM_paper_arxiv_v4}.  

\qedhere 

}

We are now ready to prove Lemma~\ref{corollary:SBMthreshold:onlyaboutdegrees}.

\proof{[of Lemma~\ref{corollary:SBMthreshold:onlyaboutdegrees}]

Let $\alpha> \beta$ be constants satisfying condition~\eqref{eq:conditionsalphabetaSBMcor_no}. Given $\Delta > 0$, we want to show that, with high probability
\begin{equation}\label{eq:tobeshownhighprobSBM5}
 \min_{i}\left( \deg_{in}(i) - \deg_{out}(i) \right) \geq \frac{\Delta}{\sqrt{\log n}} \frac{n}2(p-q). 
\end{equation}

Let us fix $i$ throughout the rest of the proof. It is clear that we can write
\[
 \deg_{in}(i) - \deg_{out}(i) = \left(\sum_{i=1}^{\frac n2-1}W_i\right) - \left(\sum_{i=1}^{n/2}Z_i\right) = \sum_{i=1}^{n/2} \left(  W_i - Z_i  \right) + Z_{\frac{n}2},
\]
where $W_1,\dots, W_m$ are i.i.d.\ $\mathrm{Bernoulli}(p)$ and $Z_1,\dots, Z_m$ are i.i.d.\ $\mathrm{Bernoulli}(q)$, independent of $W_1,\dots, W_m$. Hence, since
\[
\frac{\Delta}{\sqrt{\log n}} \left( \frac{n}2(p-q)\right) = \Delta \sqrt{\log n} \left(\frac{\alpha - \beta}2 \right),
\]
the probability of $ \deg_{in}(i) - \deg_{out}(i) <  \frac{\Delta}{\sqrt{\log n}} \left( \frac{n}2(p-q)\right) $ is equal to
\[
\PP\left[   \sum_{i=1}^{n/2} \left(  Z_i - W_i  \right) - Z_{\frac{n}2}     >  - \Delta \sqrt{\log n} \left(\frac{\alpha - \beta}2 \right)   \right] 
\]
which is upper bounded by,
\[
\PP\left[   \sum_{i=1}^{n/2} \left(  Z_i - W_i  \right)    >  - \Delta \sqrt{\log n} \left(\frac{\alpha - \beta}2 \right)  \right].
\]
Take $\Delta' = \Delta \left(\frac{\alpha - \beta}2\right) + 1$ and recall Definition~\ref{def:definitionofT}, then
\begin{eqnarray*}
 \PP\left[   \deg_{in}(i) - \deg_{out}(i) <  \frac{\Delta}{\sqrt{\log n}}  \frac{n}2(p-q)     \right]
\\ \leq  T\left(\frac{n}2,\frac{\alpha\log n}{n},\frac{\beta\log n}{n},-\Delta' \sqrt{\log n}\right)
\\ \leq  \exp\left[ -\left(\frac{\alpha + \beta}2 - \sqrt{\alpha\beta} - \delta(n)\right)\log n\right],
\end{eqnarray*}
where $\lim_{n \to infty}\delta(n)=0$, and the last inequality used Lemma~\ref{lemma:controlingTforSBM}.

Via a simple union bound, it is easy to see that,
\begin{eqnarray*}
 \PP\left[   \min_{i}\left( \deg_{in}(i) - \deg_{out}(i) \right) < \frac{\Delta}{\sqrt{\log n}} \frac{n}2(p-q)    \right] \\ 
 \leq \exp\left[ -\left(\frac{\alpha + \beta}2 - \sqrt{\alpha\beta} - 1 - \delta(n)\right)\log n\right],
\end{eqnarray*}
which means that, as long as $\frac{\alpha + \beta}2 - \sqrt{\alpha\beta}>1$, \eqref{eq:tobeshownhighprobSBM5} holds with high probability. Straightforward algebraic manipulations show that \eqref{eq:conditionsalphabetaSBMcor} implies this condition, concluding the proof of the Corollary.

\qedhere 

}

\section{Proof of the main result}\label{section:proofs}

We will prove Theorems~\ref{theorem:main} and~\ref{theorem:2} through a few Lemmas. Let us define $X$ as the non-diagonal part of $-L$ and $y\in\RR^n$ as $y = \diag\left(D_X\right)$, meaning that $y=\diag(L)$. Then $L = D_X - X$. We will separately lower bound $\max_i{y_i}$ and upper bound $\|X\|$. The upper bound on $\|X\|$ is obtained by a direct application of a result in~\cite{Bandeira_NARandomMatrixBound}.

\begin{lemma}[Remark 3.13 in~\cite{Bandeira_NARandomMatrixBound}]
\label{corollary:Master}
Let $X$ be the $n\times n$ symmetric matrix with independent centered entries.  Then there exists a universal constant $c'$, such that for every $t\ge 0$
\begin{equation}\label{eq:Master}
	\mathbf{P}[\|X\|> 3\sigma +t]
	\le ne^{-t^2/c'\sigma_{\infty}^2},
\end{equation}
where we have defined
$$
	\sigma := \max_i
	\sqrt{\sum_j \EE[X_{ij}^2]},
	\qquad\quad
	\sigma_{\infty} := \max_{ij}\|X_{ij}\|_{\infty}.
$$
\end{lemma}

Before continuing with the proof let us recall the main idea: Lemma~\ref{corollary:Master} gives that, with high probability,
\[
 \|X\| \lesssim \sigma + \sigma_{\infty}\sqrt{\log n},
\]
where $X$ is the off-diagonal part of $-L$. One the other hand, $L_{ii} = \sum_{j\in [n]\setminus i} X_{ij}$ has variance $\sigma^2$. The Central Limit Theorem would thus suggest that $L_{ii}$ behave like a gaussian of variance $\sigma^2$. Since different sums only share a single summand they are ``almost'' independent which by itself would suggest that $\max_{i}L_{ii} \sim \sigma\sqrt{\log n}$, which would imply the theorems. The proof that follows makes this argument precise.

We turn our attention to a lower bound on $\max_i{y_i}$. Recall that $y_i = \sum_{j=1}^n X_{ij}$. More specifically, we are looking for an upper bound on
\[
 \PP\left[ \max_i y_i < t \right],
\]
for a suitable value of $t$. We note that, if the $y_i$'s were independent then this could be easily done via lower bounds on the upper tail of each $y_i$. Furthermore, if the random variable $y_i$ were gaussian, obtaining such lower bounds would be trivial. Unfortunately, the random variables in question are neither independent nor gaussian, forcing major adaptations to this argument. In fact, we will actually start by lower bounding

\[
 \EE \max_{i\in [n]} y_i.
\]

We will obtain such a bound via a comparison (using Jensen's inequality) with the maximum among certain independent random variables.

\begin{lemma}
 Let $\III$ and $\JJJ$ be disjoint subsets of $[n]$. For $i\in \III$ define $z_i$ as
\begin{equation}\label{eq:defzi}
 z_i = \sum_{j\in\JJJ} X_{ij}.
\end{equation}
Then
\[
 \EE \max_{i\in [n]}y_i \geq \EE \max_{i\in \III}z_i.
\]
\end{lemma}

\proof{
\[
\EE \max_{i\in [n]}y_i = \EE \max_{i\in [n]}\sum_{j=1}^n X_{ij} \geq \EE \max_{i\in \III} \sum_{j=1}^n X_{ij}.
\]

Since $\III\cap\JJJ=\emptyset$, $\{X_{ij}\}_{i\in\III,j\in\JJJ}$ is independent from $\{X_{ij}\}_{i\in\III,j\notin\JJJ}$, and so Jensen's inequality gives
\[
 \EE \max_{i\in \III} \sum_{j=1}^n X_{ij} \geq \EE \max_{i\in \III} \left[ \sum_{j\in \JJJ} X_{ij} + \sum_{j\notin \JJJ} \EE X_{ij} \right] = \EE \max_{i\in \III} \sum_{j\in \JJJ} X_{ij} = \EE \max_{i\in \III} z_i.
\]

\qedhere

}

The following Lemma guarantees the existence of sets $\III$ and $\JJJ$ with desired properties.

\begin{lemma}\label{lemma:constructionofzi}
 There exist $\III$ and $\JJJ$ disjoint subsets of $[n]$ such that
\[
 |\III| \geq \frac18 n,
\]
and, for every $i\in \III$,
\[
 \EE z_i^2 \geq  \frac18 \sigma^2 ,
\]
where $z_i$ is defined, as in \eqref{eq:defzi}, to be $z_i = \sum_{j\in\JJJ} X_{ij}$.
\end{lemma}

\proof{

Given the matrix $X$, we start by constructing a weighted graph on $n$ nodes such that $w_{ij} = \EE X_{ij}^2$ (note that $w_{ii}=0$, for al $i$). Let $(S,S^c)$ be a partition of the vertices of this graph, with $|S|\geq \frac{n}2$, that maximizes the cut 
\[
\sum_{i\in S,\, j\in S^c}w_{ij}.
\]
It is easy to see that the maximum cut needs to be at least half of the total edge weights\footnote{One can build such a cut by consecutively selecting memberships for each node in a greedy fashion as to maximize the number of incident edges cut, see~\cite{Sahni_Gonzalez_76}.}. This readily implies
\[
\sum_{i\in S,\, j\in S^c}w_{ij} \geq \frac12 \sum_{i<j}w_{ij} = \frac14 \sum_{i\in [n]} \sum_{j\in [n]}w_{ij} = \frac14 \sum_{i\in [n]} \sum_{j\in [n]}\EE X_{ij}^2 = \frac14 n\sigma^2.
\]

Consider $z_i$, for $i\in S$, defined as
\[
 z_i = \sum_{j\in S^c}X_{ij}.
\]

We proceed by claiming that the set $\III \subset S$ of indices $i\in S$ for which
\[
 \EE z_i^2 \geq  \frac1{8} \sigma^2 ,
\]
satisfies $|\III| \geq \frac18 n$. Thus, taking $\JJJ = S^c$ would establish the Lemma.

To justify the claim, note that 
\[
  \sum_{i \in S}\EE z_i^2 = \sum_{i\in S,\, j\in S^c}w_{ij} \geq  \frac14 n\sigma^2,
\]
and
\[
  \sum_{i \in S}\EE z_i^2 \leq |\III| \max_{i\in S}\EE z_i^2 + \left(|S|-|\III|\right)  \frac1{8} \sigma^2 \leq \left(|\III|+ \frac1{8}|S|\right)\sigma^2 \leq \left(|\III|+ \frac1{8}n\right)\sigma^2,
\]
implying that $\left(|\III|+ \frac1{8}n\right)\sigma^2 \geq \frac14 n\sigma^2$.

\qedhere

}

We now proceed by obtaining a lower bound for $\EE \max_{i\in \III} z_i$, where $\III$ and $z_i$ are defined to satisfy the conditions in Lemma \ref{lemma:constructionofzi}. We note that at this point the random variables $z_i$ are independent and each is a sum of independent random variables. We use Lemma 8.1 of \cite{LT91} (for a fixed constant $\gamma = 1$) to obtain a lower bound on the upper tail of each $z_i$.

\begin{lemma}\label{lemma81LT91}\emph{[Lemma 8.1 of \cite{LT91}]}
In the setting described above, there exist two universal positive constants $K$ and $\eps$ such that for every $t$ satisfying $t\geq K \frac{\sigma}8 $ and $t \leq \eps\frac{\sigma^2}{\sqrt{8}\sigma_\infty}$, we have (for every $i\in\III$ separately)
\[
\PP \left[ z_i > t \right] \geq \exp\left( -8\frac{t^2}{\sigma^2}  \right).
\]
\end{lemma}

We are now ready to establish a lower bound on $\EE \max_{i\in[n]} y_i$.

\begin{lemma}\label{lemma:LBEmaxyi}
In the setting described above, there exist two universal positive constants $K$ and $\eps$ such that for every $t$ satisfying $t\geq K \frac{\sigma}8 $ and $t \leq \eps\frac{\sigma^2}{\sqrt{8}\sigma_\infty}$, we have
\[
 \EE \max_{i\in[n]} y_i \geq t - \left(t+n\sigma_{\infty}\right)\exp\left( -\frac{n}{\exp\left( \frac{8t^2}{\sigma^2} \right)} \right)
\]
\end{lemma}

\proof{

Let $K$ and $\eps$ be the universal constants in Lemma~\ref{lemma81LT91} and $t$ such that $K \frac{\sigma}8 \leq t \leq \eps\frac{\sigma^2}{\sqrt{8}\sigma_\infty}$. Lemma~\ref{lemma81LT91} guarantees that, for any $i\in\III$,
\[
 \PP \left[ z_i > t \right] \geq \exp\left( -8\frac{t^2}{\sigma^2}  \right).
\]

Due to the independence of the random variables $z_i$, we have
\begin{eqnarray*}
\PP \left[ \max_{i\in\III} z_i \leq t \right] & =  & \prod_{i\in \III} \PP \left[ z_i \leq t \right]  =   \prod_{i\in \III} \left( 1 - \PP \left[ z_i > t \right] \right) \\
& \leq & \left( 1 - \frac1{\exp\left( 8\frac{t^2}{\sigma^2}\right)} \right)^{|\III|} \leq  \left( 1 - \frac1{\exp\left( 8\frac{t^2}{\sigma^2}\right)} \right)^{n/8}\\
& \leq & \exp\left(-\frac{n/8}{\exp\left(8\frac{t^2}{\sigma^2}\right)}\right)
\end{eqnarray*}
where the second to last inequality follows from the fact that $|\III|\geq \frac18 n$ and the last from the fact that $\left(1-\frac1x\right)^x\leq \exp(-1)$ for $x>1$.

Since $\|X_{ij}\|_{\infty} \leq \sigma_{\infty}$ we have that, almost surely, $z_i\geq -(n-1)\sigma_{\infty}$. Thus,
\[
\EE \max_{i\in[n]} y_i \geq \EE \max_{i\in\III} z_i \geq t\left[ 1 - \exp\left(-\frac{n/8}{\exp\left(8\frac{t^2}{\sigma^2}\right)}\right) \right] - (n-1)\sigma_{\infty}\exp\left(-\frac{n/8}{\exp\left(8\frac{t^2}{\sigma^2}\right)}\right),
\]
which establishes the Lemma. 

\qedhere

}

The last ingredient we need is a concentration result to control the lower tail of $\max_{i\in[n]} y_i$ by controling its fluctuations around $\EE \max_{i\in[n]} y_i$. We make use of a result in~\cite{Massart_constants}.

\begin{lemma}\label{lemma:fromMassart}
In the setting described above, define $v$ as
\begin{equation}\label{eq:def:v}
 v = \EE\left[\max_{i\in [n]} \sum_{j=1}^n \left(X_{ij}-X'_{ij}\right)^2  \right],
\end{equation}
where $X'$ is an independent identically distributed copy of $X$.

Then, for any $x>0$:
\[
 \PP\left[ \max_{i\in[n]} y_i \leq  \EE\left[ \max_{i\in[n]} y_i \right] - x    \right] \leq \exp\left(-\frac{x^2}{7(v+\sigma_\infty x)}\right).
\]

\end{lemma}

\proof{ This Lemma is a direct consequence of Theorem 12 in~\cite{Massart_constants} by taking the independent random variables to be $Y_{(i,j)}$ such that $Y_{(i,j),t}=X_{ij}$ if $t=i$ and $Y_{(i,j),t}=0$ otherwise.  We note that there is a small typo (in the definition of the quantity $v$) in the Theorem as stated in~\cite{Massart_constants}.
\qedhere

}

At this point we need an upper bound on the quantity $v$ defined in~\eqref{eq:def:v}. This is the purpose of the following Lemma.

\begin{lemma}\label{lemma:UBforv}
In the setting above, let $X'$ is an independent identically distributed copy of $X$, then
\[ 
\EE\left[\max_{i\in [n]} \sum_{j=1}^n \left(X_{ij}-X'_{ij}\right)^2  \right] \leq 9 \sigma^2 + 90\sigma_\infty^2\log n.
\]
\end{lemma}

\proof{

We apply a Rosenthal-type inequality from Theorem~8 of~\cite{Boucheron_MomentIneq}, for each $i\in [n]$ separately, and get, for any integer $p$ and $0<\delta<1$,
\begin{eqnarray}
 \left\| \sum_{j=1}^n \left(X_{ij}-X'_{ij}\right)^2   \right\|_{p } & \leq & (1+\delta)\EE\left[ \sum_{j=1}^n \left(X_{ij}-X'_{ij}\right)^2 \right] + \frac{2p}\delta \left\| \max_{j\in [n]}  \left(X_{ij}-X'_{ij}\right)^2  \right\|_{p} \nonumber \\
 & \leq & 2(1+\delta)\sigma^2 + \frac{8p}\delta \sigma_\infty^2.
\end{eqnarray}

It is easy to see that
\[
 \EE\left[\max_{i\in [n]} \sum_{j=1}^n \left(X_{ij}-X'_{ij}\right)^2  \right] \leq n^{\frac1p} \left\| \sum_{j=1}^n \left(X_{ij}-X'_{ij}\right)^2   \right\|_{p }.
\]

Thus, taking $p=\lceil \alpha \log n \rceil$ for some $\alpha >0$ gives
\begin{eqnarray*}
  \EE\left[\max_{i\in [n]} \sum_{j=1}^n \left(X_{ij}-X'_{ij}\right)^2  \right] & \leq & n^{\frac1{\lceil \alpha \log n \rceil}}2(1+\delta)\sigma^2 + n^{\frac1{\lceil \alpha \log n \rceil}}\frac{8\lceil \alpha \log n \rceil}\delta \sigma_\infty^2 \\
& \leq &e^{\frac1\alpha}2(1+\delta)\sigma^2 + e^{\frac1\alpha}\frac{8\lceil \alpha \log n \rceil}\delta \sigma_\infty^2.
\end{eqnarray*}

Taking, for example, $\delta = 0.5$ and $\alpha = 1$ gives
\[
 \EE\left[\max_{i\in [n]} \sum_{j=1}^n \left(X_{ij}-X'_{ij}\right)^2  \right]  \leq 9 \sigma^2 + 90 \sigma_\infty^2 \log n.
\]

\qedhere

}

We now collect all our bounds in a master Lemma.

\begin{lemma}\label{lemma:mastertx}
In the setting described above, there exist universal constants $K>0$ and $\eps>0$ such that, for any $t$ satisfying $ K \frac{\sigma}8 \leq t \leq \eps\frac{\sigma^2}{\sqrt{8}\sigma_\infty} $, we have

\[
 \PP\left[ \max_{i\in[n]} y_i \leq  \frac{t}2 - \left(t+n\sigma_{\infty}\right)\exp\left( \frac{-n}{\exp\left( \frac{8t^2}{\sigma^2} \right)} \right)   \right] 
 \leq \exp\left(\frac{-t^2/10^4}{\sigma^2 + \sigma_\infty^2\log n+\sigma_\infty t}\right)
\]

\end{lemma}

\proof{

Let $t>0$ satisfy the hypothesis of the Lemma, and $x>0$.
 
Recall that Lemma~\ref{lemma:fromMassart} gives
\[
 \PP\left[ \max_{i\in[n]} y_i \leq  \EE\left[ \max_{i\in[n]} y_i \right] - x    \right] \leq \exp\left(-\frac{x^2}{7(v+\sigma_\infty x)}\right).
\]
On the other hand, Lemma~\ref{lemma:LBEmaxyi} and~\ref{lemma:UBforv} control, respectively, $\EE\left[ \max_{i\in[n]} y_i \right] $ and $v$, giving
\[
  \EE\left[ \max_{i\in[n]} y_i \right] \geq t - \left(t+n\sigma_{\infty}\right)\exp\left( -\frac{n}{\exp\left( \frac{8t^2}{\sigma^2} \right)} \right),
\]
and

\[ 
v \leq 9 \sigma^2 + 90\sigma_\infty^2\log n.
\]

Combining all these bounds,
\begin{eqnarray*}
 \PP\left[ \max_{i\in[n]} y_i \leq  t - \left(t+n\sigma_{\infty}\right)\exp\left( -\frac{n}{\exp\left( \frac{8t^2}{\sigma^2} \right)} \right) - x    \right]  \\
 \leq \exp\left(-\frac{x^2}{7(9 \sigma^2 + 90\sigma_\infty^2\log n+\sigma_\infty x)}\right).
\end{eqnarray*}

Taking $x=t/2$ establishs the Lemma.

\qedhere

}

At this point, the proofs of Theorems~\ref{theorem:main} and~\ref{theorem:2} will consist essentially of applying Lemma~\ref{lemma:mastertx} for appropriate values of $t$.

\proof{[of Theorem~\ref{theorem:main}]

Let $\beta>0$ be a constant to be defined later. Taking $t = \beta \sigma \sqrt{\log n}$ in Lemma~\ref{lemma:mastertx} gives that, in the setting described above,
\begin{eqnarray*}
 \PP\left[ \max_{i\in[n]} y_i \leq  \frac{\beta}2 \sigma \sqrt{\log n} - \left(\beta \sigma \sqrt{\log n}+n\sigma_{\infty}\right)\exp\left( -n^{1-8\beta^2} \right)   \right] 
 \\ \leq \exp\left(\frac{-\beta^2 \sigma^2 \log n/10^4}{\sigma^2 + \sigma_\infty^2\log n+\sigma_\infty (\beta \sigma \sqrt{\log n})}\right)
 \\ =  \exp\left(\frac{-\beta^2 \log n/10^4}{1 + \left(\frac{\sigma_\infty}{\sigma}\right)^2\log n+\frac{\sigma_\infty}{\sigma} \beta \sqrt{\log n}}\right)
\\ =  n^{ - \left( \frac{\beta^2/10^4}{1 + \left(\frac{\sigma_\infty}{\sigma}\right)^2\log n+\frac{\sigma_\infty}{\sigma} \beta \sqrt{\log n}} \right) },
\end{eqnarray*}
provided that  $K \frac{\sigma}8 \leq \beta \sigma \sqrt{\log n} \leq \eps\frac{\sigma^2}{\sqrt{8}\sigma_\infty}$, where $K$ and $\eps$ are the universal constants in Lemma~\ref{lemma:mastertx}.

We start by noting that, if $0 <\beta < \frac1{\sqrt{8}}$ independent of $n$, then, for $n$ large enough (not depending on $\sigma$ or $\sigma_{\infty}$),
\[
 \left(\beta \sigma \sqrt{\log n}+n\sigma_{\infty}\right)\exp\left( -n^{1-8\beta^2} \right)  \leq \frac\beta6\sigma\sqrt{\log n}.
\]

Thus, provided that  $\frac{K}{8\sqrt{\log n}} \leq \beta \leq \min\left\{ \eps\frac{\sigma}{\sqrt{8\log n}\sigma_\infty}, \frac13 \right\}$,
\[
\PP\left[ \max_{i\in[n]} y_i \leq  \frac{\beta}3 \sigma \sqrt{\log n}  \right] \leq
 n^{ - \left( \frac{\beta^2/10^4}{1 + \left(\frac{\sigma_\infty}{\sigma}\right)^2\log n+\frac{\sigma_\infty}{\sigma} \beta \sqrt{\log n}} \right) }.
\]

Let $c$ be the constant in the hypothesis of the theorem, then $\sigma > c\sqrt{\log n}\sigma_\infty$.

Let $\beta = \min\left\{ \frac{\eps c}{\sqrt{8}}, \frac13 \right\}$. Clearly, for $n$ large enough,
\[
\frac{K}{8\sqrt{\log n}} \leq \min\left\{ \frac{\eps c}{\sqrt{8}}, \frac13 \right\} \leq \min\left\{ \eps\frac{\sigma}{\sqrt{8\log n}\sigma_\infty}, \frac13 \right\},
\]
and
\[
\PP\left[ \max_{i\in[n]} y_i \leq  \min\left\{ \frac{\eps c}{6\sqrt{2}}, \frac19 \right\}  \sigma \sqrt{\log n}  \right] \leq 
 n^{ - \left( \frac{10^{-4}}{\max\left\{ \frac{8}{\eps^2 c^2}, 9 \right\} + \max\left\{ \frac{8}{\eps^2}, 9c^2 \right\}+\max\left\{ \frac{\sqrt{8}}{\eps}, 3c \right\} } \right) }.
\]

This implies that there exist constants $c_1', C_1'$ and $\beta_1'$ such that
\[
\PP\left[ \max_{i\in[n]} L_{ii} \leq  C_1'  \sigma \sqrt{\log n}  \right] \leq 
 c_1'n^{ -\beta_1' }.
\]

Recall that Corollary~\ref{corollary:Master} ensures that, for a universal constant $c'$, and for every $u\ge 0$, by taking $ t = u\sigma$,
\begin{equation}\label{eq:useofmastercoru}
\mathbf{P}[\|X\|> (3+u)\sigma] \leq ne^{-u^2\sigma^2/c'\sigma_{\infty}^2}.
\end{equation}
It is easy to see that $ne^{-u^2\sigma^2/c'\sigma_{\infty}^2} \leq ne^{-u^2(\log n)c/c'} = n^{1-u^2c/c'}$. Taking $u = \sqrt{2c'/c}$ gives
\[
\mathbf{P}\left[\|X\|> \left(3+\sqrt{2c'/c}\right)\sigma\right] \leq n^{-1}.
\]

This means that, with probability at least $1 - c_1'n^{ -\beta_1' } -  n^{-1}$ we have
\[
 \|X\| < \left(3+\sqrt{2c'/c}\right)\sigma \leq \frac{3+\sqrt{2c'/c}}{C_1' \sqrt{\log n}}\max_{i\in [n]}L_{ii},
\]
which, together with the fact that $\lambda_{\max}(L) \leq \|X\| + \max_{i\in [n]}L_{ii}$, establishes the theorem.

\qedhere

}

\proof{[of Theorem~\ref{theorem:2}]

If $\sigma > \sqrt{\log n} \sigma_\infty$ then the result follows immediately from Theorem~\ref{theorem:main}. For that reason we restrict our attention to the instances with $\sigma \leq  \sqrt{\log n} \sigma_\infty$. We start by setting 
\begin{equation}\label{eq:deftproofthm2}
t = 2\sigma \left( \frac{\sigma}{\sigma_\infty} \right)^{\frac12} (\log n)^{\frac18}.
\end{equation}
Recall that there exist $c$ and $\gamma>0$ such that
\(
 \sigma \geq c\left( \log n\right)^{\frac14 + \gamma}\sigma_{\infty},
\)
or equivalently
\[
 \frac{\sigma}{\sigma_\infty} \geq c\left( \log n\right)^{\frac14 + \gamma}.
\]

This guarantees that, for $n$ large enough (not depending on $\sigma$ or $\sigma_{\infty}$), the conditions in Lemma~\ref{lemma:mastertx} are satisfied. In fact,
\[
\frac{K\sigma}8 \leq 2\sigma \sqrt{c}\left( \log n\right)^{\frac14 + \frac{\gamma}2} \leq 2\sigma \sqrt{\frac{\sigma}{\sigma_\infty}} (\log n)^{\frac18} \leq  \frac{\eps \sigma}{\sqrt{8}} \sqrt{\frac{\sigma}{\sigma_\infty}} \sqrt{c}  \left( \log n\right)^{\frac18 + \frac{\gamma}2}  \leq \frac{\eps\sigma^2}{\sqrt{8}\sigma_\infty}.
\]

Hence, Lemma~\ref{lemma:mastertx} gives, for $t$ as in \eqref{eq:deftproofthm2},

\[
 \PP\left[ \max_{i\in[n]} y_i \leq  \frac{t}2 - \left(t+n\sigma_{\infty}\right)\exp\left( \frac{-n}{\exp\left( \frac{8t^2}{\sigma^2} \right)} \right)   \right] 
 \leq \exp\left(\frac{-t^2/10^4}{\sigma^2 + \sigma_\infty^2\log n+\sigma_\infty t}\right).
\]

We proceed by noting that, for $t=2\sigma \left( \frac{\sigma}{\sigma_\infty} \right)^{\frac12} (\log n)^{\frac18}$ and $n$ large enough (not depending on $\sigma$ or $\sigma_{\infty}$),
\[
\left(t+n\sigma_{\infty}\right)\exp\left( \frac{-n}{\exp\left( \frac{8t^2}{\sigma^2} \right)} \right)  \leq \frac{t}6.
\]

In fact, since $\sigma \leq \sigma_\infty \sqrt{\log n}$,
\[
\exp\left( \frac{-n}{\exp\left( \frac{8\left( 2\sigma \left( \frac{\sigma}{\sigma_\infty} \right)^{1/2} (\log n)^{1/8} \right)^2}{\sigma^2} \right)} \right) \leq \exp\left( \frac{-n}{\exp\left( 32 (\log n)^{3/4} \right)} \right) ,
\]
decreases faster than any polynomial.

Hence, since $t \geq 2\sigma \sqrt{c}\left( \log n\right)^{\frac14 + \frac{\gamma}2}$,

\[
 \PP\left[ \max_{i\in[n]} y_i \leq  \frac23\sigma \sqrt{c}\left( \log n\right)^{\frac14 + \frac{\gamma}2}  \right] 
 \leq \exp\left(\frac{-\left( 2\sigma \left( \frac{\sigma}{\sigma_\infty} \right)^{\frac12} (\log n)^{\frac18} \right)^2/10^4}{\sigma^2 + \sigma_\infty^2\log n+\sigma_\infty 2\sigma \left( \frac{\sigma}{\sigma_\infty} \right)^{\frac12} (\log n)^{\frac18}}\right).
\]
We proceed by noting that
\[
\frac{\left( 2\sigma \left( \frac{\sigma}{\sigma_\infty} \right)^{\frac12} (\log n)^{\frac18} \right)^2/10^4}{\sigma^2 + \sigma_\infty^2\log n+\sigma_\infty 2\sigma \left( \frac{\sigma}{\sigma_\infty} \right)^{\frac12} (\log n)^{\frac18}}
  = 
\frac{4(\log n)^{\frac14}/10^4}{\frac{\sigma_\infty}{\sigma} + \left(\frac{\sigma_\infty}{\sigma}\right)^3\log n+2\left(\frac{\sigma_\infty}{\sigma}\right)^{\frac32}(\log n)^{\frac18}}
\]
Since $\frac{\sigma_\infty}{\sigma} \leq \frac1c(\log n)^{-\frac14-\gamma}$, we have that, for $n$ large enough and a constant $c''$
\[
 \PP\left[ \max_{i\in[n]} y_i \leq  \frac23\sigma \sqrt{c}\left( \log n\right)^{\frac14 + \frac{\gamma}2}  \right] 
 \leq \exp\left(-c'' (\log n)^{\gamma}\right).
\]

At this point we upper bound $\|X\|$, as in the proof of Theorem~\ref{theorem:main}. Recall, as in \eqref{eq:useofmastercoru}, for any $u>0$,
\[
\mathbf{P}[\|X\|> (3+u)\sigma] \leq ne^{-\frac{u^2\sigma^2}{c'\sigma_\infty^2}}.
\]
Hence,
\[
\mathbf{P}[\|X\|> (3+u)\sigma] \leq ne^{-\frac{u^2 c^2}{c'} (\log(n))^{\frac12 + 2\gamma}} .
\]
Taking $u = (\log n)^{\frac14}$ gives
\[
\mathbf{P}[\|X\|> \left(3+(\log n)^{\frac14}\right)\sigma] \leq e^{-\frac{c^2}{c'} (\log(n))^{2\gamma}}.
\]
The rest of the proof follows the final arguments in the proof of Theorem~\ref{theorem:main}.

\qedhere

}

\section{Conclusion and future directions}

Theorems~\ref{theorem:main} and~\ref{theorem:2} are valid for matrices whose entries may be distributed in very different ways. This potentially allows one to use them in order to obtain strong guarantees for deterministically censored versions of the problems described, where the measurements are obtained only for edges of a deterministic graph (a similar model was studied, for example, in~\cite{Abbe_Z2Synch}).

The problem of recovery in the stochastic block model with multiple balanced clusters, also referred to as multisection, is a natural generalization of the one considered here and also admits a semidefinite relaxation. While the results here do not seem to be directly applicable in the analysis of that algorithm, in part because the construction of a dual certificate in that setting is considerably more involved, some of the ideas in the present paper can be adapted for the estimates needed there. These also provide interpretable, and sharp, guarantees. We refer the interested reader to~\cite{Agawarl_multisection_SBM}.

Regarding directions for future investigations, from the random matrix side of things it would be interesting to investigate what happens when $\sigma \gg \sigma_{\infty}$ but $\frac{\sigma}{\sigma_{\infty}} = o\left((\log n)^{\frac14}\right)$, as this setting is not captured by our results. It would be particularly interesting also to understand whether analogues of these results exist for instances where the off-diagonal entries of $L$ are not independent~\footnote{For the particular example of connectivity of an Erd\H{o}s--R\'{e}nyi graph, it is possible to use the matrix concentration approach~\cite{Tropp:TailBoundsRM,Tropp_randommatrix_monograph} to obtain a guarantee that, while being a factor away from optimal, appears to be adaptable to instances where edges have particular types of dependencies --- we refer the reader to Section 5.3. in the monograph~\cite{Tropp_randommatrix_monograph}.}. 

From the point of view of applications, a natural question is which other semidefinite relaxations have these optimality guarantees. A general understanding in that direction would be remarkable.

\subsection*{Acknowledgements}

The author acknowledges Amit Singer, Emmanuel Abbe, Ramon van Handel, and Georgina Hall for many insightful discussions on the topic of this paper and, specialy, for motivating the author
 to write this manuscript. Many thanks to Ramon van Handel for his crucial help locating the references for the strongest version of the theorems used in the proof of our main results.

The author is also indebted to Amit Singer, Ramon van Handel, Dustin G. Mixon, Nicolas Boumal, and Joel Tropp for valuable comments on early drafts of this manuscript.

The author presented most of these results in various seminars throughout the end of 2014 and beginning of 2015. Many questions and comments raised by the audience greatly improved the quality of this manuscript, a warm thanks to all of them.

\bibliographystyle{plain}
\bibliography{../afonso.bib}

\begin{thebibliography}{10}

\bibitem{Abbe_Z2Synch}
E.~Abbe, A.~S. Bandeira, A.~Bracher, and A.~Singer.
\newblock Decoding binary node labels from censored edge measurements: Phase
  transition and efficient recovery.
\newblock {\em Network Science and Engineering, IEEE Transactions on},
  1(1):10--22, Jan 2014.

\bibitem{Abbe_Z2SynchER}
E.~Abbe, A.~S. Bandeira, A.~Bracher, and A.~Singer.
\newblock Linear inverse problems on {E}rd{\H{o}}s-{R}\'enyi graphs:
  Information-theoretic limits and efficient recovery.
\newblock {\em IEEE International Symposium on Information Theory (ISIT2014)},
  2014.

\bibitem{SBM_paper_arxiv_v4}
E.~Abbe, A.~S. Bandeira, and G.~Hall.
\newblock Exact recovery in the stochastic block model.
\newblock {\em Available online at arXiv:1405.3267v4 [cs.SI]}, 2014.

\bibitem{Agawarl_multisection_SBM}
N.~Agarwal, A.~S. Bandeira, K.~Koiliaris, and A.~Kolla.
\newblock Multisection in the stochastic block model using semidefinite
  programming.
\newblock {\em Submitted}, 2015.

\bibitem{Alizadeh_93}
F.~Alizadeh.
\newblock Interior point methods in semidefinite programming with applications
  to combinatorial optimization.
\newblock {\em SIAM Journal on Optimization}, 5:13--51, 1993.

\bibitem{Anderson_Guionnet_Zeitouni_IntroRandomMatrices}
G.~W. Anderson, A.~Guionnet, and O.~Zeitouni.
\newblock {\em An introduction to random matrices}.
\newblock Cambridge studies in advanced mathematics. Cambridge University
  Press, Cambridge, New York, Melbourne, 2010.

\bibitem{Bandeira_rankrecoveryangsynch}
A.~S. Bandeira, N.~Boumal, and A.~Singer.
\newblock Tightness of the maximum likelihood semidefinite relaxation for
  angular synchronization.
\newblock {\em Available online at arXiv:1411.3272 [math.OC]}, 2014.

\bibitem{Bandeira2014open}
A.~S. Bandeira, Y.~Khoo, and A.~Singer.
\newblock Open problem: Tightness of maximum likelihood semidefinite
  relaxations.
\newblock In {\em Proceedings of the 27th Conference on Learning Theory},
  volume~35 of {\em JMLR W\&CP}, pages 1265--1267, 2014.

\bibitem{Bandeira_Singer_Spielman_OdCheeger}
A.~S. Bandeira, A.~Singer, and D.~A. Spielman.
\newblock A {C}heeger inequality for the graph connection {L}aplacian.
\newblock {\em SIAM J. Matrix Anal. Appl.}, 34(4):1611--1630, 2013.

\bibitem{Bandeira_NARandomMatrixBound}
A.~S. Bandeira and R.~v.~Handel.
\newblock Sharp nonasymptotic bounds on the norm of random matrices with
  independent entries.
\newblock {\em Annals of Probability, to appear}, 2015.

\bibitem{correlation_clustering_04}
N.~Bansal, A.~Blum, and S.~Chawla.
\newblock Correlation clustering.
\newblock {\em Machine Learning}, 56(1-3):89--113, 2004.

\bibitem{Boppana_SBM1}
R.~B. Boppana.
\newblock Eigenvalues and graph bisection: An average-case analysis.
\newblock {\em In 28th Annual Symposium on Foundations of Computer Science},
  pages 280--285, 1987.

\bibitem{Boucheron_MomentIneq}
S.~Boucheron, O.~Bousquet, G.~Lugosi, and P.~Massart.
\newblock Moment inequalities for functions of independent random variables.
\newblock {\em Ann. Probab.}, 33(2):514--560, 2005.

\bibitem{Bryc_06_MarkovMatrices}
W.~Bryc, A.~Dembo, and T.~Jiang.
\newblock Spectral measure of large random {H}ankel, {M}arkov and {T}oeplitz
  matrices.
\newblock {\em The Annals of Probability}, 34(1):pp. 1--38, 2006.

\bibitem{Chen_Goldsmith_ISIT2014}
Y.~Chen and A.~J. Goldsmith.
\newblock Information recovery from pairwise measurements.
\newblock {\em IEEE International Symposium on Information Theory (ISIT2014)},
  2014.

\bibitem{Chen_Goldsmith_2015}
Y.~Chen, C.~Suh, and A.~J. Goldsmith.
\newblock Information recovery from pairwise measurements: A shannon-theoretic
  approach.
\newblock {\em Available online at arXiv:1504.01369 [cs.IT]}, 2015.

\bibitem{Chung_ComplexGraphsNetworks}
F.~Chung and L.~Lu.
\newblock {\em Complex Graphs and Networks (Cbms Regional Conference Series in
  Mathematics)}.
\newblock American Mathematical Society, Boston, MA, USA, 2006.

\bibitem{FanChung_SpectralGraphTheory}
F.~R.~K. Chung.
\newblock {\em Spectral Graph Theory}.
\newblock AMS, 1997.

\bibitem{Cucuringu_Z2Synch}
M.~Cucuringu.
\newblock Synchronization over {Z}2 and community detection in signed multiplex
  networks with constraints.
\newblock {\em Journal of Complex Networks}, 2015.

\bibitem{Szarek_Davidson_survey}
K.~Davidson and S.~Szarek.
\newblock Local operator theory, random matrices and {B}anach spaces.
\newblock In {\em Handbook on the Geometry of Banach spaces}, volume~1, pages
  317--366. Elsevier Science, 2001.

\bibitem{Decelle_SBM}
A.~Decelle, F.~Krzakala, C.~Moore, and L.~Zdeborov\'a.
\newblock Asymptotic analysis of the stochastic block model for modular
  networks and its algorithmic applications.
\newblock {\em Phys. Rev. E}, 84, December 2011.

\bibitem{Ding_RandomLaplacians}
X.~Ding and T.~Jiang.
\newblock Spectral distribution of adjacency and {L}aplacian matrices of random
  graphs.
\newblock {\em The Annals of Applied Probability}, 20(6):2086--2117, 2010.

\bibitem{Durrett_RandomGraphDynamics}
R.~Durrett.
\newblock {\em Random Graph Dynamics (Cambridge Series in Statistical and
  Probabilistic Mathematics)}.
\newblock Cambridge University Press, New York, NY, USA, 2006.

\bibitem{ErdosRenyi_59}
P.~Erd\H{o}s and A.~R\'{e}nyi.
\newblock {On random graphs, I}.
\newblock {\em Publicationes Mathematicae (Debrecen)}, 6:290--297, 1959.

\bibitem{Feige_Kilian_bisection_01}
U.~Feige and J.~Kilian.
\newblock Heuristics for semirandom graph problems.
\newblock {\em Journal of Computer and System Sciences}, 63(4):639 -- 671,
  2001.

\bibitem{MXGoemans_DPWilliamson_1995}
M.~X. Goemans and D.~P. Williamson.
\newblock Improved approximation algorithms for maximum cut and satisfiability
  problems using semidefine programming.
\newblock {\em Journal of the Association for Computing Machinery},
  42:1115--1145, 1995.

\bibitem{Grotschel_Lovasz_Schrijver_ellipsoidmethod}
M.~Gr\"otschel, L.~Lov\'asz, and A.~Schrijver.
\newblock The ellipsoid method and its consequences in combinatorial
  optimization.
\newblock {\em Combinatorica}, 1(2):169--197, 1981.

\bibitem{Hajek_et_al_SBM_SDP}
B.~Hajek, Y.~Wu, and J.~Xu.
\newblock Achieving exact cluster recovery threshold via semidefinite
  programming.
\newblock {\em Available online at arXiv:1412.6156}, 2014.

\bibitem{Hajek_et_al_SBM_SDP_extensions}
B.~Hajek, Y.~Wu, and J.~Xu.
\newblock Achieving exact cluster recovery threshold via semidefinite
  programming: Extensions.
\newblock {\em Available online at arXiv:1502.07738}, 2015.

\bibitem{SKhot_2002}
S.~Khot.
\newblock On the power of unique 2-prover 1-round games.
\newblock {\em Thiry-fourth annual ACM symposium on Theory of computing}, 2002.

\bibitem{LT91}
M.~Ledoux and M.~Talagrand.
\newblock {\em Probability in {B}anach spaces}, volume~23 of {\em Ergebnisse
  der Mathematik und ihrer Grenzgebiete (3) [Results in Mathematics and Related
  Areas (3)]}.
\newblock Springer-Verlag, Berlin, 1991.

\bibitem{Lovasz_ShannonCapacity}
L.~Lovasz.
\newblock On the shannon capacity of a graph.
\newblock {\em IEEE Trans. Inf. Theor.}, 25(1):1--7, 1979.

\bibitem{Massart_constants}
P.~Massart.
\newblock About the constants in {T}alagrand's concentration inequalities for
  empirical processes.
\newblock {\em The Annals of Probability}, 28(2), 2000.

\bibitem{Massoulie_SBM}
L.~Massouli{\'e}.
\newblock Community detection thresholds and the weak ramanujan property.
\newblock In {\em Proceedings of the 46th Annual ACM Symposium on Theory of
  Computing}, STOC '14, pages 694--703, New York, NY, USA, 2014. ACM.

\bibitem{McSherry_SBM1}
F.~McSherry.
\newblock Spectral partitioning of random graphs.

\bibitem{Mossel_SBM3_exact}
E.~Mossel, J.~Neeman, and A.~Sly.
\newblock Consistency thresholds for the planted bisection model.
\newblock {\em Available online at arXiv:1407.1591v2 [math.PR]}, July 2014.

\bibitem{Mossel_SBM2}
E.~Mossel, J.~Neeman, and A.~Sly.
\newblock A proof of the block model threshold conjecture.
\newblock {\em Available online at arXiv:1311.4115 [math.PR]}, January 2014.

\bibitem{Mossel_SBM1}
E.~Mossel, J.~Neeman, and A.~Sly.
\newblock Stochastic block models and reconstruction.
\newblock {\em Probability Theory and Related Fields (to appear)}, 2014.

\bibitem{Nesterov_Nemirovskii_interiorpointmethods_94}
Y.~Nesterov and A.~Nemirovskii.
\newblock {\em Interior-Point Polynomial Algorithms in Convex Programming}.
\newblock Society for Industrial and Applied Mathematics, 1994.

\bibitem{Raghavendra_2008_optimalitySDP_UG}
P.~Raghavendra.
\newblock Optimal algorithms and inapproximability results for every {CSP}?
\newblock In {\em Proceedings of the Fortieth Annual ACM Symposium on Theory of
  Computing}, STOC '08, pages 245--254. ACM, 2008.

\bibitem{Sahni_Gonzalez_76}
S.~Sahni and T.~Gonzalez.
\newblock P-complete approximation problems.
\newblock {\em J. ACM}, 23(3):555--565, July 1976.

\bibitem{ASinger_2011_angsync}
A.~Singer.
\newblock Angular synchronization by eigenvectors and semidefinite programming.
\newblock {\em Appl. Comput. Harmon. Anal.}, 30(1):20 -- 36, 2011.

\bibitem{Tao_topicsRMT}
T.~Tao.
\newblock {\em Topics in Random Matrix Theory}.
\newblock Graduate studies in mathematics. American Mathematical Soc., 2012.

\bibitem{Tropp:TailBoundsRM}
J.~A. Tropp.
\newblock User-friendly tail bounds for sums of random matrices.
\newblock {\em Foundations of Computational Mathematics}, 12(4):389--434, 2012.

\bibitem{Tropp_randommatrix_monograph}
J.~A. Tropp.
\newblock An introduction to matrix concentration inequalities.
\newblock {\em Found. Trends Mach. Learning}, 8(1--2):1--230, 2015.

\bibitem{LVanderberghe_SBoyd_1996}
L.~Vanderberghe and S.~Boyd.
\newblock Semidefinite programming.
\newblock {\em SIAM Review}, 38:49--95, 1996.

\bibitem{VershyninNARandomMatrices}
R.~Vershynin.
\newblock Introduction to the non-asymptotic analysis of random matrices.
\newblock {\em Chapter 5 of: Compressed Sensing, Theory and Applications.
  Edited by Y. Eldar and G. Kutyniok. Cambridge University Press}, 2012.

\bibitem{Wigner_1958}
E.~P. Wigner.
\newblock On the distribution of the roots of certain symmetric matrices.
\newblock {\em Annals of Mathematics}, 67(2):pp. 325--327, 1958.

\end{thebibliography}

\end{document}